\documentclass[a4paper]{article}

\usepackage{graphicx}
\usepackage{amsmath}
\usepackage{amsthm}
\usepackage{amssymb}
\usepackage{amsfonts}
\usepackage{latexsym}
\usepackage{epsfig}
\usepackage{color}

\newcommand{\comment}[1]{}

\comment{
\textwidth 7.0 in
\textheight 9.0 in
\oddsidemargin -0.5 in
\evensidemargin -0.5 in
\headheight -1.0in
\footskip 0.2 in
\parskip .08in
\parindent 0in
}

\newcommand{\beql}[1]{\begin{equation}\label{#1}}
\newcommand{\eeq}{\end{equation}}
\newcommand{\Ds}{\displaystyle}

\newcommand{\Abs}[1]{{\left|{#1}\right|}}
\newcommand{\Lone}[1]{{\left\|{#1}\right\|_1}}

\newcommand{\Linf}[1]{{\left\|{#1}\right\|_\infty}}
\newcommand{\Norm}[1]{{\left\|{#1}\right\|}}

\newcommand{\Floor}[1]{{\left\lfloor{#1}\right\rfloor}}

\newcommand{\Qed}{\hfill \mbox{\raggedright \rule{.07in}{.1in}}}

\newcommand{\Set}[1]{{\left\{{#1}\right\}}}

\newcommand{\RR}{{\mathbb R}}
\newcommand{\CC}{{\mathbb C}}
\newcommand{\ZZ}{{\mathbb Z}}
\newcommand{\NN}{{\mathbb N}}
\newcommand{\TT}{{\mathbb T}}

\newcommand{\One}[1]{{\bf 1}\left(#1\right)}

\newcommand{\inner}[2]{{\langle #1, #2 \rangle}}
\newcommand{\Inner}[2]{{\left\langle #1, #2 \right\rangle}}
\newcommand{\dens}{{\rm dens\,}}

\newcommand{\supp}{{\rm supp\,}}
\newcommand{\dist}{{\rm dist\,}}
\newcommand{\vol}{{\rm vol\,}}
\newcommand{\ft}[1]{\widehat{#1}}
\newcommand{\FT}[1]{\left(#1\right)^\wedge}
\newcommand{\nozero}[1]{{#1\setminus\Set{0}}}

\newcounter{dfn}
\def\thedfn{\arabic{dfn}}

\newcounter{obs}
\def\theobs{\arabic{obs}}

\newcounter{thm}
\newcounter{mysec}
\newcounter{mysubsec}[mysec]
\newtheorem{theorem}{Theorem}[section]
\newtheorem{conjecture}{Conjecture}[section]

\newtheorem{lemma}{Lemma}[section]

\newtheorem{definition}{Definition}[section]
\newtheorem{remark}{Remark}[section]
\newtheorem{open-problem}{Open Problem}

\numberwithin{equation}{section}

\begin{document}

\begin{center}
{\bf\Large The study of translational tiling with Fourier Analysis}\\
\ \\
{\large Lectures given at the Workshop on Fourier Analysis and Convexity,}\\
{\large Universit\`a di Milano--Bicocca},\\
{\large June 11-22, 2001}\\
\vskip 3cm
{\sc Mihail N. Kolountzakis\footnote{
Supported in part by European Commission IHP Network HARP
(Harmonic Analysis and Related Problems),
Contract Number: HPRN-CT-2001-00273 - HARP.}}\\
\ \\
Department of Mathematics,\\
University of Crete,\\
Knossos Ave,\\
714 09 Iraklio,\\
Greece.\\
\ \\
E-mail: {\tt kolount@member.ams.org}\\
\ \\
WWW: {\tt http://fourier.math.uoc.gr/\~\,mk}\\
\vskip 3cm
March 2003
\end{center}

\pagebreak

\tableofcontents


\pagebreak
\noindent{\bf Forward}
\addcontentsline{toc}{section}{Forward}

In this survey I will try to describe how Fourier Analysis is used in the
study of translational tiling.
Right away I will emphasize two restrictions
that separate this area from the general theory of tilings.

\begin{figure}[t] \centering
\input{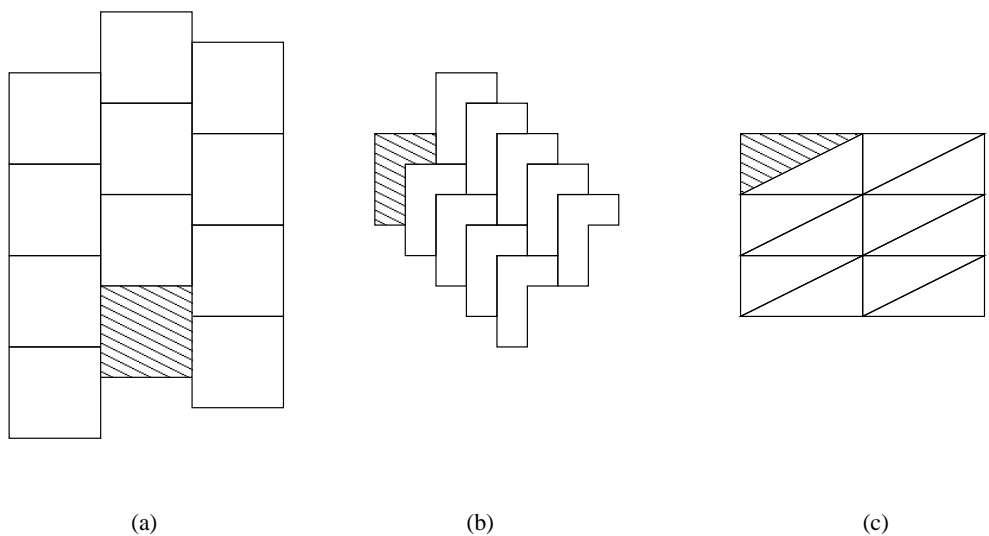}
\caption{Examples of tiling with the shaded objects. In (c)
a tiling by a triangle is shown that is using rotations as well as translations.
We will not deal with such tilings here. In (a) a tiling by a square
is shown and in (b) a tiling by an L-shaped region. In (b) the set
of translations is a lattice, but not in (a).
}
\ \\ \hrule
\label{fig:basic-examples}
\end{figure}

\begin{itemize}
\item 
There is only one tile.
This is an object that is moved around
in space (whatever space we are trying to tile, most generally an abelian
group) in a way that there
are no ``overlaps'' among the several copies of it and
almost nothing, in the sense of Lebesgue or counting measure, is left uncovered.
This object may be a
domain in space or a function defined on space, usually nonnegative.
Examples are shown in Figure \ref{fig:basic-examples}.
\item
The only allowed motions of the tile are translations.
No rotations or reflections of the object are allowed.
In fancier language, we are tiling abelian
groups, not vector spaces.
\end{itemize}
This paper is broken up into three ``lectures'', which correspond roughly
to the three hour-long lectures I gave in the Universit\'a di Milano--Bicocca,
in June 2001, during the meeting on {\em Fourier Analysis and Convexity}.
Lecture 1 has to do with how Fourier Analysis is used to prove
structure, or rigidity, in tilings. In Lecture 2, some problems are presented
about lattice-tiling and in Lecture 3 a tiling problem of Functional Analysis
is discussed, the Fuglede Conjecture on spectral domains.

An advance apology: I will describe mostly material with which I am aquainted
the most, through my own work.

Finally, I would like to thank the organizers L. Brandolini, L. Colzani, 
A. Iosevich and G. Travaglini for organizing this great meeting and giving
me the chance to participate.

\pagebreak
\section{Lecture 1: Introduction to the method and structure of tilings.}

\subsection{Tiling and density}

It's time for the first definition, of what tiling means. We speak mostly
of tiling $\RR^d$ and $\ZZ^d$ in this paper, but tiling makes sense on
all abeliean groups.
\begin{figure}[h] \centering
\input{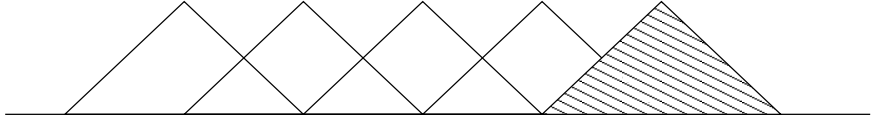}
\caption{A triangle function tiling the real line}
\ \\ \hrule
\label{fig:soft-tile}
\end{figure}
\begin{definition}
{\bf (Translational tiling)}\\
Suppose $0\le f \in L^1(\RR^d)$ and $\Lambda\subseteq\RR^d$ is a discrete multiset.
We say that $f$ tiles $\RR^d$ with $\Lambda$ at level (or weight) $\ell$ if
$$
\sum_{\lambda\in\Lambda} f(x-\lambda) = \ell,  \ \ \mbox{a.e.($x$)}.
$$
We write: $f+\Lambda = \ell\,\RR^d$.
\end{definition}
In Figure \ref{fig:soft-tile} a tiling by the triangle function
$f(x) = (1-\Abs{x})^+$ is shown with translation set $\Lambda = \ZZ$
and level $1$.
In the particular case when $f = \chi_\Omega$ is the indicator function
of a measurable domain $\Omega \subseteq \RR^d$ of finite measure, we
write also $\Omega + \Lambda = m \RR^d$, where the positive integer $m$
represents the level of the (generally multiple) tiling.

The tiling assumption $f+\Lambda = \ell\RR^d$ has some immediate
implications about the density properties of the multiset $\Lambda$.
\begin{definition}\label{def:density}
{\bf (Density)}\\
A multiset $\Lambda\subseteq\RR^d$ has asymptotic density
$\rho$ if
$$
\lim_{R\to\infty} {\#(\Lambda \cap B_R(x)) \over \Abs{B_R(x)}} \to
\rho
$$
uniformly in $x\in\RR^d$. We write $\rho = \dens \Lambda$.\\

We say that $\Lambda$ has (uniformly) bounded density if the fraction
above is bounded by a constant $\rho$ uniformly for $x\in\RR$ and
$R>1$.
We say then that $\Lambda$ has density (uniformly) bounded by $\rho$.

Last, the upper density of a set $\Lambda\subseteq\RR^d$ is defined as
$$
\limsup_{R\to\infty} \sup_{x\in\RR^d} {\#(\Lambda \cap B_R(x)) \over \Abs{B_R(x)}}.
$$
\end{definition}
\begin{remark}
According to this definition a set $\Lambda$ may have {\em density
uniformly bounded by a number $\rho<\infty$} yet $\dens \Lambda$ may not exist.
\end{remark}

\begin{lemma}\label{lm:tiling-implies-bounded-density}
If $0\le f \in L^1(\RR^d)$ is not the zero function and $f+\Lambda = \ell\RR^d$
then $\Lambda$ has bounded density.
\end{lemma}
\noindent{\bf Proof.}
By hypothesis
$$
\sum_{a \in \Lambda} f(x-a) = \ell,\ \ \ \mbox{almost everywhere},
$$
and clearly $\ell > 0$.
Choose $R > 1$ so that $J = \int_{B_R(0)} f > 0$, where $B_R(0)$ is the
ball centered at $0$ with radius $R$.
Let $t\in\RR^d$ be arbitrary. We have
\begin{eqnarray*}
\Abs{B_{2R}(0)} \cdot \ell & = & \int_{B_{2R}(t)} \sum_{a \in \Lambda} f(x-a)~ dx \\
 & \geq & \int_{B_{2R}(t)} \sum_{\Abs{a-t} < R} f(x-a)~ dx \\
 & \geq & \#(\Lambda \cap B_R(t)) \int_{B_R(0)} f .
\end{eqnarray*}
Thus $\#(\Lambda \cap B_R(t)) \leq \Abs{B_{2R}(0)}\ell/J$ is bounded independent of $t$,
which implies that $\Lambda$ has uniformly bounded density.
\Qed

Working similarly on easily gets the following lemma.
\begin{lemma}\label{lm:tiling-implies-density}
If $0\le f \in L^1(\RR^d)$ is not the zero function and $f+\Lambda = \ell\RR^d$
then $\Lambda$ density $\dens \Lambda = \ell (\int f)^{-1}$.
\end{lemma}

It is time also to define {\em packing}.
\begin{definition}
{\bf (Packing)}\\
Suppose $0\le f \in L^1(\RR^d)$ and $\Lambda\subseteq\RR^d$ is a discrete multiset.
We say that $f$ packs $\RR^d$ with $\Lambda$ at level $\ell$ if
$$
\sum_{\lambda\in\Lambda} f(x-\lambda) \le \ell,  \ \ \mbox{a.e.($x$)}.
$$
We write: $f+\Lambda \le \ell\,\RR^d$.
\end{definition}

The following lemma is almost trivial, yet useful.
\begin{lemma}\label{lm:packing-implies-density-bound}
If $0\le f \in L^1(\RR^d)$ is not the zero function and $f+\Lambda \le \ell\RR^d$
is a packing then $\Lambda$ has density uniformly bounded by $\ell (\int f)^{-1}$.
\end{lemma}
Finally, one can easily prove the following about translation sets.
\begin{lemma}\label{lm:tiling-sets-are-separated}
Suppose $f+\Lambda \le \ell \RR^d$ and ${\rm esssup}\ f = \ell$.
Then
\beql{separation}
\inf\Set{\Abs{\lambda-\mu}:\ \lambda,\mu \in \Lambda, \lambda \neq \mu} > 0.
\eeq
In particular, if $E+\Lambda = \RR^d$ is a tiling by the set $E$ at level 1
then \eqref{separation} holds.
\end{lemma}

\subsection{Tiling in Fourier space}\label{sec:tiling-in-fourier-space}

\begin{figure}[h] \centering
\input{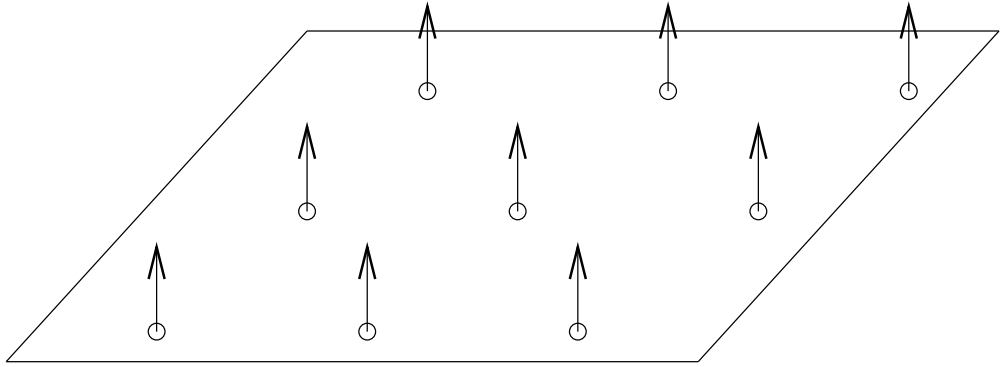}
\caption{The measure $\delta_\Lambda$ corresponding to some $\Lambda$ in the plane.}
\ \\ \hrule
\label{fig:delta-lambda}
\end{figure}

Next, we associate to any point multiset $\Lambda$ the measure
$$
\delta_\Lambda = \sum_{\lambda\in\Lambda} \delta_\lambda,
$$
where $\delta_\lambda$ is one unit point mass at the point $\lambda$
(see Figure \ref{fig:delta-lambda}).
Generally, this measure is infinite globally but has finite
total variation in any bounded set, at least when the set $\Lambda$
has bounded density. This is the case whenever $\Lambda$ is involved
in a tiling. It follows that
$$
\Abs{\delta_\Lambda}(B_R(t)) \le C R^d,
$$
which implies that the object $\delta_\Lambda$ is a so-called
tempered distribution, a bounded linear functional on the Schwarz
space ${\cal S}$ of smooth functions which, along with all their partial
derivatives, decay faster than any power at infinity.

If $T$ is a tempered distribution one defines its Fourier Transform
$\ft T$ by duality as follows:
$$
\ft{T}(\phi) = T(\ft\phi),
$$
for any $\phi\in{\cal S}$ (it is easy to prove that the Fourier
Transform $\ft\phi$ is also in ${\cal S}$).
We normalize the Fourier Transform for a function $f \in L^1(\RR^d)$ as
$$
\ft f(t) = \int e^{-2\pi i \inner{t}{x}}f(x)~dx,
$$ 
which leads to the inversion formula
$$
f(x) = \int e^{2\pi i \inner{t}{x}}\ft f (t)~dt,
$$
whenever $\ft f \in L^1$, which happens for all functions $f \in {\cal S}$.

We are now in the position to argue formally as follows.
Suppose $f + \Lambda = \ell\RR^d$. This means that
$$
\sum_{\lambda\in\Lambda} f(x-\lambda) = \ell, \ \ \mbox{(a.e. $x$)},
$$
which we rewrite as a convolution
$$
f * \delta_\Lambda = \ell.
$$
Take the Fourier Transform of both sides to get
$$
\ft{f}\cdot\ft{\delta_\Lambda} = \ell\delta_0.
$$
As the support of the right hand side is just $\Set{0}$ we conclude that
\beql{main}
\supp \ft{\delta_\Lambda} \subseteq \Set{0} \cup {\cal Z}(\ft{f}),
\eeq
where we denote the zero-set of the continuous function $g$ by ${\cal Z}(g)$:
$$
{\cal Z}(g) = \Set{x\in\RR^d:\ g(x) = 0}.
$$
The inclusion in \eqref{main} is the starting point of the method of
applying Fourier Analysis to translational tiling.
Whenever we have tiling, we deduce \eqref{main}.
Sometimes we may be able to get tiling from \eqref{main}, but we usually
need some extra conditions to make this conclusion.

Having argued formally, let us now prove carefully the following theorem.
Notice that we have essentially added the condition $\ft f \in C^\infty$
to make the argument go through. This condition is automatically valid whenever
$f$ has compact support, as, for instance, when $f$ is the indicator function
of a bounded domain (the classical geometric situation), but will definitely
not be there when we talk about the Fuglede problem in Lecture 3. There we will
need a different theorem of this sort,
with different assumptions (see Theorem \ref{th:supp-equiv-tiling}).
\begin{theorem}\label{th:tiling-implies-condition}
Suppose that $f\in L^1(\RR^d)$ is nonnegative, $\ft{f} \in C^\infty$
and $f+\Lambda = \ell\RR^d$
for some multiset $\Lambda$. Then \eqref{main} follows.
\end{theorem}
\noindent{\bf Proof.}
Let $K = \Set{0}\cup {\cal Z}(\ft{f})$, which is a closed set. Inclusion \eqref{main}
means (by the definition of the support of a tempered distribution) that
$\ft{\delta_\Lambda}(\psi) = 0$ for all smooth
$\psi$ supported in $K^c$ (see Figure \ref{fig:test-function}).
For such a $\psi$ 
\begin{figure}[h] \centering
\begin{center}\input{test-function.pstex_t}\end{center}
\caption{A test function $\psi$ supported away from $\Set{0}\cup\Set{\ft f = 0}$}
\ \\ \hrule
\label{fig:test-function}
\end{figure}
\begin{eqnarray}
\left( \widetilde{\ft{f}} \, \psi \right)^\wedge (\lambda) &=& 
      \int \overline{\ft{f}(-x)} \psi(x) e^{-2\pi i\lambda x}\,dx\nonumber\\
 &=&  \int \int \overline{f(y)} e^{-2\pi i yx} \psi(x)
		e^{-2\pi i\lambda x}\,dx\,dy\nonumber\\
 &=&  \int \overline{f(y)} \int \psi(x)e^{-2\pi i(\lambda+y)x}\,dx\,dy\nonumber\\
 &=&  \int \overline{f(y)} \ft{\psi}(\lambda+y)\,dy\nonumber\\
 &=&  \int \overline{f(-t)} \ft{\psi}(\lambda-t)\,dt\nonumber\\
 &=& (\widetilde{f}*\ft{\psi})(\lambda)\label{convolution},
\end{eqnarray}
where we use the notation $\widetilde{f}(x) = \overline{f(-x)}$.

We must show $\ft{\delta_\Lambda}(\psi) = 0$.
We have
$$
\ft{\delta_\Lambda}(\psi) = \ft{\delta_\Lambda}\left(
	\widetilde{\ft{f}} \cdot {\psi \over \widetilde{\ft{f}} } \right).
$$
Notice that $\widetilde{\ft f}$ and $\ft f$ have the same zeros (since $f$ is real),
so the quotient $\phi = \psi / \widetilde{\ft{f}}$ is a $C^\infty_0(K^c)$
function.
We have
\begin{eqnarray*}
\ft{\delta_\Lambda}(\psi)
 &=& \ft{\delta_\Lambda} ( \widetilde{\ft f} \,\phi ) \\
 &=& \delta_\Lambda \left( ( \widetilde{\ft f} \phi )^\wedge \right)
\ \ \mbox{(by the definition of the Fourier Transform for distributions)}\\
 &=& \sum_{\lambda\in\Lambda}
	 ( \widetilde{\ft f} \phi )^\wedge (\lambda)
\ \ \mbox{(by the definition of $\delta_\Lambda$)}\\
 &=& \sum_\lambda (\widetilde{f}*\ft\phi)(\lambda)
\ \ \mbox{(by \eqref{convolution})}\\
 &=& \sum_\lambda \int \widetilde{f}(\lambda-x) \ft\phi(x)\,dx\\
 &=& \int \sum_\lambda \overline{f(x-\lambda)} \ft\phi(x)\,dx\\
 &=& \ell \int \ft\phi(x)\,dx
\ \ \mbox{(since $f+\Lambda = \ell\RR^d$)}\\
 &=& \ell\phi(0)\\
 &=& 0
\ \ \mbox{(as $0 \notin \supp\phi$)}.
\end{eqnarray*}
\Qed

\subsubsection{The lattice case and the sufficiency of the support condition for tiling}
\label{sec:lattice}

Suppose $\Lambda = A\ZZ^d$, $A \in GL(d, \RR)$, is a lattice
in $\RR^d$ (a discrete subgroup which contains $d$ linearly
independent vectors).
The Fourier Transform of the tempered distribution $\delta_\Lambda$
takes a particularly simple form as claimed by the {\em Poisson Summation Formula:}
\beql{psf}
\ft{\delta_\Lambda} = {1\over\det A} \delta_{\Lambda^*},
\eeq
where
$$
\Lambda^*=\Set{\xi\in\RR^d:\ \inner{\xi}{\lambda} \in\ZZ,
\ \ \forall \lambda\in\Lambda} = A^{-\top}\RR^d
$$
is the {\em dual lattice} of $\Lambda$ (see Figure \ref{fig:dirac-combs}).
\begin{figure}[h] \centering
\begin{center}\input{dirac-combs.pstex_t}\end{center}
\caption{The ``Dirac comb'' $\delta_\Lambda$ when $\Lambda = {1\over2}\ZZ$,
and its Fourier Transform, the comb $2\delta_{2\ZZ}$.}
\ \\ \hrule
\label{fig:dirac-combs}
\end{figure}

The Poisson Summation Formula is usually stated as the 
equality
$$
\sum_{\lambda\in\Lambda}\ft\phi(\lambda) =
 {1\over \det A} \sum_{\lambda^*\in\Lambda^*} \phi(\lambda^*),
$$
for all $\phi \in {\cal S}$, and this is exactly the content of
\eqref{psf}, as the Fourier Transform of $\delta_\Lambda$ is defined by duality.

Equation \eqref{main} now gives the implication below, valid
for any lattice $\Lambda$,
$$
\mbox{$f+\Lambda$ is a tiling} \Rightarrow
 \mbox{$\ft f$ vanishes on $\Lambda^*\setminus\Set{0}$}.
$$
This is in fact easy to prove using ordinary multiple Fourier Series,
after applying a linear transformation that maps $\Lambda$ to $\ZZ^d$.
Working this way one gets easily that the above implication is, in fact,
an equivalence, so that
\beql{lattice-tiling-in-fourier-space}
\mbox{$f+\Lambda$ is a tiling} \Leftrightarrow
 \mbox{$\ft f$ vanishes on $\Lambda^*\setminus\Set{0}$}.
\eeq
We prefer however to stick to using \eqref{main} as our guiding tool and not
mention Fourier Series.
As to why the reverse implication holds, the answer is in the following theorem.
\begin{theorem}\label{th:converse}
Suppose that $\Lambda$ is a multiset of bounded density and that $f$ is a nonnegative
integrable function on $\RR^d$.
Suppose also that $\ft{\delta_\Lambda}$ is locally a measure
and that
$$
\supp\ft{\delta_\Lambda} \subseteq \Set{0} \cup \Set{\ft f = 0}.
$$
Then $\Lambda$ has density and $f+\Lambda = \ell\RR^d$,
for $\ell = \int f \cdot \dens\Lambda$.
\end{theorem}
Intuitively, to kill a tempered distribution which is a measure any zero (of
whatever order) suffices.

\noindent{\bf Proof.}
Let $F(x) = \sum_{\lambda\in\Lambda}f(x-\lambda)$. We want to show that $F$
is a constant $\ell$ and for this it is enough to show that for any nonnegative
$\ft\phi \in {\cal S}$ we have $\int F \ft \phi = \ell \int\ft\phi = \ell \phi(0)$.
We have
\begin{eqnarray*}
\int F \ft\phi &=& \sum_\lambda \int f(x-\lambda)\ft\phi(x)~dx\\
 &=& \int f(y) \sum_\lambda \ft\phi(y+\lambda)~dy\\
 &=& \int f(y) \delta_\Lambda(\ft\phi(y-\cdot))~dy\\
 &=& \int f(y) \ft{\delta_\Lambda} ( e^{2\pi i yx} \phi(x) )~dy\\
 &=& \int\int f(y) e^{2\pi i yx} \phi(x)~d\ft{\delta_\Lambda}(x)~dy\\
 &=& \int \ft f(-x) \phi(x) ~d\ft{\delta_\Lambda}(x)\\
 &=& \ft{\delta_\Lambda}(\Set{0})\phi(0)\ft f(0),
\end{eqnarray*}
which proves the desired equality with $\ell = \int f \cdot \ft{\delta_\Lambda}(\Set{0})$.
The fact that $\Lambda$ has density and the value for
$\dens \Lambda$ follow from Lemma \ref{lm:tiling-implies-density}.
\Qed

\subsection{Structure of tilings in dimension 1}

We can now show the following theorem \cite{leptin-mueller,kolountzakis-lagarias}.
\begin{theorem}\label{th:structure-in-dim-one}
{\bf (Leptin and M\" uller, 1991, and Kolountzakis and Lagarias, 1996)}\\
Suppose $0\le f\in L^1(\RR)$ and has compact support.
Suppose also that
$$
f + \Lambda = \ell \RR,
$$
for some $\Lambda \in \RR$.
Then there are $J \in \NN$, $\alpha_j, \beta_j \in \RR$, $j=1,\ldots,J$,
$\alpha_j > 0$, such that
$$
\Lambda = \bigcup_{j=1}^J (\alpha_j\RR+\beta_j).
$$
\end{theorem}
That is, tiling sets for compactly supported tiles in dimension $1$ are finite
unions of complete arithmetic progressions.

\begin{remark}
Only recently (after the first draft of this survey was written in March 2003)
has the author learned of the work of Leptin and M\" uller \cite{leptin-mueller}.
In \cite{leptin-mueller} Theorem \ref{th:structure-in-dim-one} was proved,
by essentially the same method.
\end{remark}

\subsubsection{The idempotent theorem, the Bohr group and Meyer's theorem}

This extreme structure is, in the end, a consequence of
P.J. Cohen's idempotent theorem on a general abelian group \cite{cohen}.
\begin{theorem}\label{th:cohen}
{\bf (Cohen, 1959)}\\
If $\mu \in M(G)$ is a finite measure on a locally compact abelian
group $G$, such that $\ft\mu$ takes only finitely many values then, for
any such value $c$,
the set $S = \Set{\gamma\in\ft G:\ \ft\mu(\gamma)=c}$ belongs to the
open coset ring of $\ft G$.
\end{theorem}
The (open) coset ring is defined below.
\begin{definition}\label{dfn:coset-ring}
{\bf (The coset ring of a group)}\\
The {\em coset ring} of an abelian group $G$ is the smallest
collection of subsets of $G$ which is closed under finite unions,
finite intersections and complements
and which contains all cosets of $G$.
For a topological group $G$ the smallest ring
of subsets of $G$ which contains all {\em open}
cosets is called the {\em open coset ring} of $G$.
\end{definition}
Cohen's theorem therefore says that $S$
can be constructed with finitely many set-theoretic
operations from the open cosets of $\ft G$.

The group $\ft G$ is called the {\em dual group} of $G$ and
is the group of continuous characters on $G$, that is, the group
of all group homomorphisms $G \to \CC$ with the group operation
beeing the pointwise multiplication.
It can be proved that $\ft{\ft G}$ is isomorphic (as a topological
group) with $G$ (Pontryagin duality) and that $\ft G$ is compact
if and only if $G$ is discrete.
Further $\ft{G\times H} = \ft G \times \ft H$.
Some dual group pairs are the following: $(\ZZ,\TT)$, $(\RR,\RR)$,
$(\ZZ_n,\ZZ_n)$, $(\RR^d,\RR^d)$, $(\ZZ^d,\TT^d)$.

If $\mu$ is a finite measure on $G$ its Fourier Transform is a continuous
function on $\ft G$ defined by
$$
\ft\mu(\xi) = \int_G \overline{\xi(x)}~d\mu(x),
$$
the integration carried out with respect to the essentially unique
translation invariant measure on $G$ called the {\em Haar measure}.
For example, when $G=\RR$ the Haar measure is Lebesgue measure and
$\xi(x) = e^{2\pi I \xi x}$.
(The reader should consult \cite{rudin:groups} for the basic definitions
and facts about Fourier Analysis on locally compact abelian groups.)

We do not use Cohen's theorem directly, but rather a consequence of
it discovered by Y. Meyer \cite{meyer}. 
\begin{theorem}\label{th:meyer}
{\bf (Meyer, 1970)}\\
Let $\Lambda \subseteq \RR^d$ be a discrete set and $\delta_\Lambda$
be the Radon measure
$$
\delta_\Lambda = \sum_{\lambda \in \Lambda} c_\lambda \delta_\lambda,
\ \ c_\lambda\in S,
$$
where $S\subseteq\nozero{\CC}$ is a finite set.
Suppose that $\delta_\Lambda$ is tempered,
and that $\ft{\delta_\Lambda}$ is a Radon measure on $\RR^d$ which
satisfies
\beql{growth}
\Abs{\ft{\delta_\Lambda}} (B_R(0)) \le C_1 R^d, ~ \mbox{as} ~ R
\rightarrow
  \infty,
\eeq
where $C_1>0$ is a constant.
Then, for each $s\in S$,
the set
$$
\Lambda_s = \Set{\lambda\in\Lambda:\ c_\lambda = s}
$$
is in the coset ring of $\RR^d$.
\end{theorem}

\noindent{\bf Proof.}
Let $\phi \in C_c^\infty (B_1(0))$, $\phi (0) = 1$,
so that its Fourier Transform satisfies
$\Abs{\ft{\phi} ( \xi )} \le C_\alpha \Abs{\xi}^{-\alpha}$ for all
$\alpha > 0$.
For positive integers $n$ define the {\em functions}
$$
\mu_n (x) = \phi (nx) * \mu(x).
$$
Their Fourier Transforms satisfy
$$
\ft{\mu}_n(\xi) = {1 \over n^d} \ft{\phi}(\xi/n) \ft{\mu}(\xi),
$$
hence the $\ft{\mu}_n$ are all measures.
We claim that the measures $\ft{\mu_n}$ are uniformly bounded measures,
i.e. $\Abs{\ft{\mu_n}} ( \RR^d ) \le C$, where $C$ is independent of $n$.
Indeed
\beql{bound-near-zero}
\Abs{\ft{\mu_n}} ( B_n(0) ) \leq {1 \over n^d} ||\ft{\phi}||_\infty\,
	\Abs{\ft{\mu}} ( B_n(0) ) \le C_1 ||\ft{\phi}||_\infty,
\eeq
by our assumption on the growth of $\Abs{\ft{\mu}} ( B_n(0) )$.

Furthermore, if $2^k \gg n$ we have
(using the fact that $\Abs{\ft{\phi} ( \xi )} \le C \Abs{\xi}^{-d-1} $
as $\xi\to\infty$)
\begin{eqnarray*}
\Abs{\ft{\mu_n}} (B_{2^{k+1}}(0) \setminus B_{2^k}(0)) & \le &
C {1 \over n^d} \Abs{\ft\phi}_{B_{2^{k+1}/n}(0) \setminus B_{2^k/n}(0)}
   \Abs{\ft\mu}(B_{2^{k+1}}(0))\\
 &\le& C {1 \over n^d} \left({2^k \over n}\right)^{-d-1} 2^{(k+1)d}\\
 &\le& C n 2^{-k}.
\end{eqnarray*}
Hence
$$
\Abs{\ft{\mu_n}}(B_n(0)^c) \le \sum_{n \le 2^k} \Abs{\ft{\mu_n}}(B_{2^{k+1}}(0) \setminus B_{2^k}(0))
 \le C n \sum_{n \le 2^k} 2^{-k} \le C 1,
$$
which, together with \eqref{bound-near-zero}, shows that the sequence
$\Abs{\ft{\mu_n}}(\RR^d)$ is bounded.

Notice also that $\lim_{n \rightarrow \infty} \mu_n (x) = c_x$ if
$x\in \Lambda$ and is 0 otherwise.
This is a consequence of the fact that $\Lambda$ is discrete
and the support of $\phi (nx)$ shrinks to 0.

We now use the following properties of $\overline{\RR^d}$,
the Bohr compactification of $\RR^d$, a locally compact abelian group.
\begin{enumerate}
\item
$\overline{\RR^d}$ is the dual group of $\RR_\delta^d$,
the $d$-dimensional Euclidean space with the discrete topology.
Therefore $\overline{\RR^d}$ is a compact group being the dual group of a discrete group.
\item
$\RR^d \subseteq \overline{\RR^d}$ as topological spaces and $\RR^d$
is dense in $\overline{\RR^d}$.
Identifying the continuous functions on $\overline{\RR^d}$
with bounded continuous functions on $\RR^d$ we get that
$$
C( \overline{\RR^d} ) \subseteq C( \RR^d ) \cap L^\infty ( \RR^d )
$$
is a Banach space inclusion.
\end{enumerate}
Since the measures $\widehat{\mu_n}$ are uniformly bounded
they act on all bounded continuous functions on $\RR^d$, and
consequently also on all continuous functions on $\overline{\RR^d}$.
That is they constitute a uniformly bounded family of linear functionals on
$C( \overline{\RR^d} )$.
By the Banach-Alaoglu theorem there exists a measure $\nu$ on $\overline{\RR^d}$
such that for every $f\in C(\overline{\RR^d})$ there is a subsequence
of $\ft{\mu_n}$, call it again $\widehat{\mu_n}$, such that
$$
\ft{\mu_n}(f) \to \nu(f), \ \mbox{ as $n\to\infty$}.
$$
Applying this with each character of $\overline{\RR^d}$ in place of $f$
we obtain that
$$
\ft{\nu} (x) = \lim_{n \to\infty} \ft{\ft{\mu_n}} (x) = c_{-x}, \ \ \mbox{if $-x\in \Lambda$},
$$
and is 0 otherwise.
Hence $\ft\nu$ has the finite range $-S$.
By Theorem \ref{th:cohen} the set $-\Lambda$,
and thus $\Lambda$, belongs to the open coset ring of $\RR_\delta^d$.
Since $\RR_\delta^d$ has the discrete topology the open coset ring
is the same as the coset ring of $\RR^d$.
\Qed

Since we need to know what kind of sets the elements of the
coset ring of $\RR^d$ are, we use the following general
theorem \cite{kolountzakis:polygons},
which says that discrete elements of the coset ring
can always be constructed from {\em discrete} cosets using
finitely many unions, intersections and complementations.
\begin{theorem}\label{th:coset-ring}
{\bf (Kolountzakis, 2000)}\\
Let $G$ be a topological abelian group and let ${\cal R}$ be
the least ring of sets which contains the discrete cosets of
$G$.
Then ${\cal R}$ contains all discrete elements of the coset ring of
$G$.
\end{theorem}
In dimension $1$ this implies the following result by Rosenthal \cite{rosenthal}.
\begin{theorem}\label{th:rosenthal}
{\bf (Rosenthal, 1966)}\\
The elements of the coset ring of $\RR$ which are discrete in the
usual topology of $\RR$ are precisely the sets of the form
\beql{ap-unions}
F \bigtriangleup \bigcup_{j=1}^J ( \alpha_j \ZZ + \beta_j ) ~,
\eeq
where $F \subseteq \RR$ is finite, $J \in \NN$, $\alpha_j > 0$ and $\beta_j \in
\RR$
($\bigtriangleup$ denotes symmetric difference).
\end{theorem}

\subsubsection{Getting structure in dimension $1$}
\label{sec:structure-in-dim-1}

In this section we prove Theorem \ref{th:structure-in-dim-one}.
Assume that $\Lambda \subset \RR$ is set of bounded density and
that $f+\Lambda = \ell \RR$ for a function $f \in L^1$
of compact support, contained in, say, $(-A,A)$.
We will use \eqref{main}, so the first thing
to do is to obtain information on the set ${\cal Z}(\ft f) = \Set{\ft f = 0}$.

We look at the Fourier Transform of $f$ defined on the complex
numbers
$$
\ft f(z) = \int_{\RR} e^{-2\pi i zx} f(x)~dx,\ \ \ (z \in \CC).
$$
Since $f$ is supported in $(-A,A)$ it follows that $\ft f$ is entire
so that ${\cal Z}(\ft f)$ is a discrete subset of $\RR$.
Furthermore $\ft f$
satisfies the growth bound
$$
\Abs{\ft{f} (z)} \leq \int_{-A}^A e^{2\pi x Im (z)}\Abs{f(x)}~dx
 \leq ||f||_1 e^{2\pi A \Abs{z}}.
$$
If $N(T)$ counts the number of zeros of $\widehat{f}(z)$ in the
disk $\Set{z :\ \Abs{z} \leq T}$, an application of Jensen's formula gives
$$
\limsup_{T \rightarrow \infty} {N(T) \over T} \leq C A.
$$
Write $B$ for the discrete set $\Set{0} \cup {\cal Z}(\ft f)$, so that by \eqref{main}
the tempered distribution $\ft \delta_\Lambda$ is supported on $B$.
It is well known, and easy to prove, that a tempered distribution
supported at a single point $b$ is necessarily a finite linear combination
of derivatives of $\delta_b$, and the same proof gives that
$$
\ft{\delta_\Lambda} = \sum_{b\in B} P_b(\partial) \delta_b.
$$
Here $P_b(\partial) = \sum_{j=0}^N c_j {\partial^j \over \partial x^j}$
is differential polynomial operator applied
on the Dirac point mass at $b$. (The degree $N$ can be taken the same for all $b \in B$
as any tempered distribution has finite degree. This is not used below.)

\noindent
{\bf Step 1} All $P_b$ are constants (hence $\ft{\delta_\Lambda}$ is locally a measure)

\noindent
Focus on a single $b\in B$ and let $\phi$ be a smooth function
of compact support.
\begin{figure}[h] \centering
\input{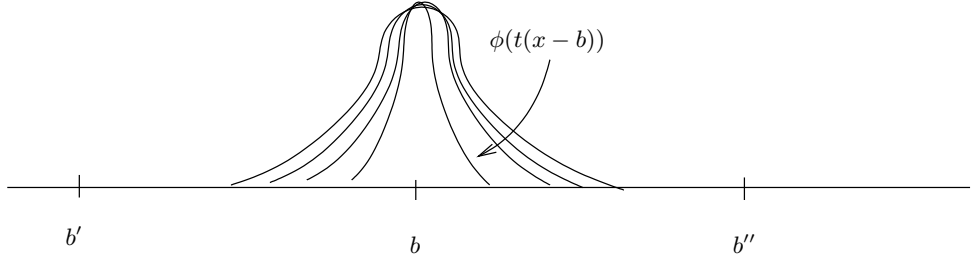}
\caption{Picking out the distribution $\ft{\delta_\Lambda}$ at $b$ by applying
it on $\phi(t(x-b))$. For large $t$ the other points of
set $B$ are left out and the behavior at $b$ is isolated.}
\ \\ \hrule
\label{fig:lens}
\end{figure}
Examine the quantity
$$
I(t) = \ft{\delta_\Lambda}\left(\phi(t(x-b))\right),\ \ \ (t\to\infty),
$$
as shown in Figure \ref{fig:lens}.
For large $t$ this equals
\begin{eqnarray*}
(P_b(\partial)\delta_b) \left(\phi(t(x-b))\right) &=&
  \left( \sum_{j=0}^N c_j \delta_b^{(j)} \right) \left(\phi(t(x-b))\right)\\
 &=& \sum_{j=1}^N c_j (-1)^j \phi^{(j)}(0) t^j.
\end{eqnarray*}
Choose $\phi^{(j)}(0) = (-1)^j$ to get the above expression equal to
$$
\sum_{j=1}^N c_j t^j.
$$
Next we will bound the growth of $I(t)$.

\noindent
Let
$$
g(x) = \phi(t(x-b)),\ \ \ 
 \ft g(\xi) = {1\over t} e^{-2\pi i b\xi/t}\ft\phi\left({\xi\over t}\right).
$$
By duality
\begin{eqnarray*}
\Abs{I(t)}
 &=& \Abs{\ft{\delta_\Lambda}(g)}\\
 &=& \Abs{\delta_\Lambda(\ft g)}\\
 &\le& {1\over t}\sum_\lambda \Abs{\ft\phi\left({\lambda\over t}\right)}\\
 &=& {1\over t}\sum_{\Abs{\lambda}\le t} + {1\over t}\sum_{\Abs{\lambda}> t}\\
 &\le& C + C \sqrt{t} \sum_{n=\Floor{t}}^\infty n^{-3/2}\\
 &=& O(\sqrt{t}).
\end{eqnarray*}
We used the bounded density of $\Lambda$ for the convergence of the sum $\sum_{\Abs{\lambda}>t}$,
and the fact that
\beql{decay}
\Abs{\ft\phi(\xi)} = O\left(\Abs{\xi}^{-M}\right)
\eeq
for any $M>0$ we wish. We took $M = 3/2$.

Since $I(t)$ cannot even grow linearly it follows that the degree $N$
is zero and we can now write
$$
\ft{\delta_\Lambda} = \sum_{b\in B} c_b \delta_b,
$$
for some constants $c_b$.

\noindent
{\bf Step 2} The coefficients $c_b$ are uniformly bounded.

\noindent
To prove this we are just
a bit more careful in the last estimate and now use a
$\phi$ which is $1$ at $0$.
For large $t$ then
$$
c_b = \ft{\delta_\Lambda}(\phi(t(x-b))),
$$
and one can get a bound for this by duality which
does not involve $t$ at all using the exponent
$M=2$ instead of $M=3/2$ in \eqref{decay}.

\noindent
{\bf Step 3} Use of Meyer's Theorem

\noindent
Now the crucial condition
$$
\Abs{\ft{\delta_\Lambda}}(-R, R) \le C R
$$
in Meyer's Theorem holds (remember there is a linear
number of zeros and at each one we have a bounded mass),
hence, by Rosenthal's Theorem \ref{th:rosenthal},
$$
\Lambda = \bigcup_{j=1}^J(\alpha_j\ZZ+\beta_j) \,\triangle\, F
$$
for some real numbers $\alpha_j, \beta_j$ and
finite set $F$.

\noindent
{\bf Step 4} $F$ is empty

\noindent
Otherwise $\ft{\delta_\Lambda}$ would have a continuous
part, a trigonometric polynomial due to $F$.
But it cannot have such a continuous part
as its support is discrete.

\begin{open-problem}
Is the main theorem true if $f$ is only supposed to be in $L^1$
but not of compact support?

What if $f$ is an indicator function?
\end{open-problem}

\subsection{Structure of some polygonal tilings in dimension $2$}

The one-dimensional tiling problem treated in the previous section
is very particular. One cannot expect this rigid structure in higher
dimension. For example, even when the tile is a square in two dimensions,
one cannot expect every tiling of it to be fully-periodic, in the sense
of posessing a period lattice of full-rank. One can, after all, make vertical
columns of squares which can be shifted vertically, within themselves, arbitrarily,
preserving the tiling property (see Figure \ref{fig:basic-examples} (a)).
It is clear that there is no horizontal period here, in general.
One might suspect that there is always, no matter what the tile, at least
one period, but this phenomenon, if true, must happen only in dimension two.
In dimension three one can construct cube tilings with no periods at all.
First make horizontal layers of cubes some of which have no period along the $x$-axis and some
others having no period along the $y$-axis. Consider these tiled slabs as rigid bodies and move
each of them by an arbitrary horizontal vector thereby destroying all vertical periods as well.

\begin{open-problem}
If $E\subseteq\RR^2$, is it true
that in any tiling $E+\Lambda = \RR^2$ the set $\Lambda$ must posess at
least one period-vector?
\end{open-problem}

The main difficulty in dimension two and higher is that the zero set
of $\ft f$ is not a discrete set any more, at least under no set of 
reasonable assumptions about $f$ (such as compact support was in dimension one).
Therefore, from our basic condition \eqref{main}
on obtains that $\ft{\delta_\Lambda}$ is supported, in general, on a subset of
the plane, which, under some reasonable assumptions, is a collection
of submanifolds of codimension one.
The structure of such distributions is much richer of course than those
supported at points, and this is the main source of difficulty, at least compared
with the one-dimensional problem.

In this section we will show the following result \cite{kolountzakis:polygons} in two-dimensions.
\begin{theorem}\label{th:convex-polygons}
{\bf (Kolountzakis, 2000)}\\
Suppose that $P$ is a symmetric convex polygon in the plane
which tiles (multiply) with the multiset $\Lambda$:
$$
P + \Lambda = m\RR^d
$$
at some integer level $m$.
If $P$ is not a parallelogram then
$\Lambda$ is a finite union of two-dimensional lattices.
\end{theorem}
The convexity assumption here is only used to guarantee
that each edge-direction appears in the polygon exactly twice.
For a more general theorem see \cite{kolountzakis:polygons}.

If one tries to use \eqref{main} directly, one encounters the
problems mentioned above, mainly the fact that the zero set ${\cal Z}(\ft{\chi_P})$
is not discrete, but rather a one-dimensional set.

\begin{figure}[h] \centering
\input{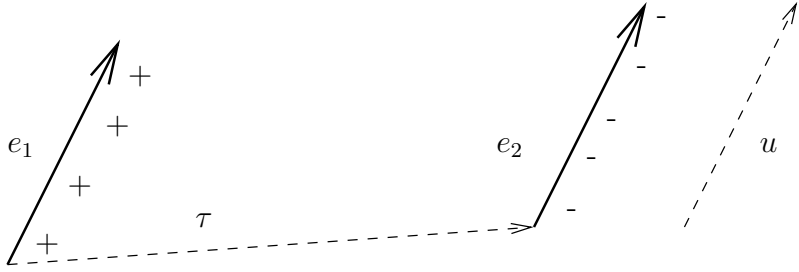}
\caption{The measure $\mu_u$ supported on two parallel edges of
the polygon $e_1$ and $e_2$, with opposite sign on each edge.}
\ \\ \hrule
\label{fig:edge-measure}
\end{figure}

Let $e_1$ and $e_2$ be two edges of the polygon $P$
of the same direction $u$. By the symmetry of $P$ they have the same length.
We can then write (here $e_1$ and $e_2$ are viewed as point-sets
in $\RR^2$ and $\tau$ as a vector)
$$
e_2 = e_1 + \tau,
$$
for some $\tau\in\RR^2$.
(For each set $A$ and vector $x$ we write $A+x = \Set{a+x:\ a\in A}$.)
Let then $\mu_u$ be the measure which is equal to arc-length on
$e_1$ and negative arc-length on $e_2$ (see Figure \ref{fig:edge-measure}).
Since every part of a translate of $e_1$ in the tiling $P+\Lambda$
has to be cancelled by part of a copy of $e_2$ it follows that
$$
\sum_{\lambda\in\Lambda} \mu_u ( x -\lambda)
$$
is the zero measure in $\RR^2$.
It also intuitively obvious that the vanishing of the above measure
for all relevant directions (i.e.\ those appearing as edge-directions) $u$
also implies tiling at some integer level.

So a convex symmetric polygon $P$ tiles multiply with a multiset $\Lambda$ if and
only if for each pair $e$ and $e+\tau$ of parallel edges of $P$
\beql{vanishes}
\sum_{\lambda\in\Lambda} \mu_e (x -\lambda) = 0,
\eeq
where $\mu_e$ is the measure in $\RR^2$ that is arc-length on
$e$ and negative arc-length on $e+\tau$.
Condition \eqref{vanishes} then becomes
$\mu_e * \delta_\Lambda = 0$ or, taking Fourier Transforms
(arguing as in \S \ref{sec:tiling-in-fourier-space}),
$$
\ft{\mu_e}\cdot\ft{\delta_\Lambda} = 0.
$$
and
\beql{containment-in-zero-set}
\supp\ft{\delta_\Lambda} \subseteq {\cal Z}(\ft{\mu_e})
\eeq
for {\em all} edge-directions $e$.

\subsubsection{The shape of the zero-set}

Here we study the zero-set of $\ft{\mu_e}$ and determine its structure.
We first calculate $\ft{\mu_e}$ in the particular
case when $e$ is parallel to the $x$-axis, for simplicity.
Let $\mu \in M(\RR^2)$ be the measure defined by duality by
$$
\mu(\phi) = \int_{-1/2}^{1/2} \phi(x, 0)~dx,
\ \ \forall \phi\in C(\RR^2).
$$
That is, $\mu$ is arc-length on the line segment joining the points
$(-1/2, 0)$ and $(1/2, 0)$. Calculation gives
$$
\ft{\mu}(\xi, \eta) = {\sin{\pi\xi} \over \pi\xi}.
$$
Notice that $\ft{\mu}(\xi, \eta) = 0$ is equivalent to
$\xi \in \nozero{\ZZ}$.

If $\mu_L$ is the arc-length measure on the line segment joining
$(-L/2,0)$ and $(L/2,0)$ we have
$$
\ft{\mu_L}(\xi, \eta) = {\sin{\pi L \xi} \over \pi\xi}
$$
and
$$
{\cal Z}(\ft{\mu_L}) = \Set{(\xi, \eta):\ \xi\in \nozero{L^{-1}\ZZ} }.
$$
Write $\tau = (a, b)$ and let $\mu_{L,\tau}$ be the measure
which is arc-length on the segment joining $(-L/2,0)$ and $(L/2,0)$
translated by $\tau/2$ and negative arc-length on the same segment
translated by $-\tau/2$.
That is, we have
$$
\mu_{L, \tau} = \mu_L * (\delta_{\tau/2} - \delta_{-\tau/2}),
$$
and, taking Fourier Transforms, we get
$$
\ft{\mu_{L,\tau}}(\xi, \eta) = -2{\sin{\pi L \xi} \over \pi\xi}
\sin\pi(a\xi+b\eta).
$$
Define $u = {\tau\over\Abs{\tau}^2}$ and $v=(1/L, 0)$.
It follows that ($u^\perp$ is a unit vector orthogonal to $u$)
$$
{\cal Z}(\ft{\mu_{L,\tau}}) = (\ZZ u + \RR u^\perp) \cup
	(\nozero{\ZZ}v + \RR v^\perp).
$$
(Each of the two summands in the union above corresponds to each of
the factors in the formula for $\ft{\mu_{L,\tau}}$.)
This a set of straight lines of direction $u^\perp$ spaced by
$\Abs{u}$ and containing $0$ plus a similar set of lines of direction
$v^\perp$, spaced by $\Abs{v}$ and containing zero. However in the
latter
set of parallel lines the straight line through $0$ has been removed
(see Figure \ref{fig:zero-set}).
\begin{figure}[h] \centering
\input{zero-set.pstex_t}
\caption{The zero-set ${\cal Z}(\ft{\mu_{L,\tau}}) = (\ZZ u + \RR u^\perp) \cup
        (\nozero{\ZZ}v + \RR v^\perp)$, with $u = {\tau\over\Abs{\tau}^2}$
and $v=(1/L, 0)$}
\ \\ \hrule
\label{fig:zero-set}
\end{figure}
We state this as a theorem for later use, formulated in a
coordinate-free way.

\begin{definition}
{\bf (Geometric inverse of a vector)}\\
The geometric inverse of a non-zero vector $u \in \RR^2$ is the vector
$$
u^* = {u \over \Abs{u}^2}.
$$
\end{definition}

\begin{theorem}\label{th:polygon-zero-set}
Let $e$ and $e+\tau$ be two parallel line segments (translated by
$\tau$, of
magnitude and direction described by $e$, symmetric with respect to
$0$).
Let also $\mu_{e,\tau}$ be the measure which charges $e$ with its
arc-length
and $e+\tau$ with negative its arc-length.
Then
\beql{zero-shape}
{\cal Z}(\ft{\mu_{e,\tau}}) = (\ZZ \tau^* + \RR \tau^{*\perp}) \cup
	(\nozero{\ZZ} e^* + \RR e^{*\perp}).
\eeq
\end{theorem}

\subsubsection{Completion of the argument}
The {\em intersection} of all the relevant ${\cal Z}(\ft{\mu_e})$
is easily shown to be a discrete set, except when $P$ is a parallelogram.

To conclude the argument we show that the tempered distribution
$\ft{\delta_\Lambda}$ is (a) locally a measure, and (b) the point masses
of $\ft{\delta_\Lambda}$ are uniformly bounded.
This is accomplished using the following two Theorems.
\begin{theorem}\label{th:at-zero}
Suppose that $\Lambda \in \RR^d$ is a multiset with density $\rho$,
$\delta_\Lambda = \sum_{\lambda\in\Lambda} \delta_\lambda$,
and that $\ft{\delta_\Lambda}$ is a measure in a neighborhood of
$0$.
Then $\ft{\delta_\Lambda}(\Set{0}) = \rho$.
\end{theorem}

\noindent{\bf Proof.}
Take $\phi \in C^\infty$ of compact support with $\phi(0) = 1$.
We have
\begin{eqnarray*}
\ft{\delta_\Lambda}(\Set{0}) &=&
 \lim_{t\to\infty} \ft{\delta_\Lambda}(\phi(t x))\\
&=& \lim_{t\to\infty} \delta_\Lambda(t^{-d}\ft\phi(\xi/t))\\
&=& \lim_{t\to\infty} t^{-d}
\sum_{\lambda\in\Lambda}\ft\phi(\lambda/t)\\
&=& \lim_{t\to\infty} \sum_{n\in\ZZ^d} \sum_{\lambda\in Q_n} t^{-d}
\ft\phi(\lambda/t)
\end{eqnarray*}
where, for fixed and large $T>0$,
$$
Q_n = [0,T)^d + T n,\ \ \ \ n \in \ZZ^d.
$$
Since $\Lambda$ has density $\rho$ it follows that for each
$\epsilon>0$
we can choose $T$ large enough so that for all $n$
$$
\Abs{\Lambda \cap Q_n} = \rho \Abs{Q_n}(1 + \delta_n),
$$
with $\Abs{\delta_n} \le \epsilon$.
For each $n$ and $\lambda\in Q_n$ we have
$$
\ft\phi(\lambda / t) = \ft\phi(T n / t) + r_\lambda
$$
with $\Abs{r_\lambda} \le C T t^{-1} 
 \Norm{\nabla\ft\phi}_{L^\infty(t^{-1}Q_n)}$.
Hence
\begin{eqnarray*}
\ft{\delta_\Lambda}(\Set{0}) &=&
 \lim_{t\to\infty} \sum_{n\in\ZZ^d} t^{-d} \sum_{\lambda\in Q_n}
   (\ft{\phi}(Tn/t) + r_\lambda) \\
&=& \lim_{t\to\infty} \sum_{n\in\ZZ^d} t^{-d} 
	\rho\Abs{Q_n}(1+\delta_n)\ft\phi(Tn/t) + \\
&\ & \lim_{t\to\infty} \sum_{n\in\ZZ^d} 
		t^{-d} \sum_{\lambda\in Q_n} r_\lambda \\
&=& \lim_{t\to\infty}S_1 + \lim_{t\to\infty}S_2.
\end{eqnarray*}
We have
\beql{tmp2}
\Abs{S_1 - \sum_n t^{-d} \rho \Abs{Q_n} \ft\phi(Tn/t)} \le
  \epsilon \sum_n t^{-d} \rho \Abs{Q_n} \Abs{\ft\phi(Tn/t)}
\eeq
The first sum in \eqref{tmp2} is a Riemann sum for 
$\rho\int_{\RR^d}\ft\phi = \rho$ and the second is
a Riemann sum for $\rho\int_{\RR^d}\Abs{\ft\phi}<\infty$.

For $S_2$ we have
\begin{eqnarray*}
\Abs{S_2} &\le&
 C \sum_{n\in\ZZ^d} t^{-d} \rho \Abs{Q_n} (1+\delta_n) T t^{-1}
   \Norm{\nabla\ft\phi}_{L^\infty(t^{-1}Q_n)} \\
&\le& C \rho T t^{-1} \sum_{n\in\ZZ^d} t^{-d}\Abs{Q_n}
   \Norm{\nabla\ft\phi}_{L^\infty(t^{-1}Q_n)}.
\end{eqnarray*}
The sum above is a Riemann sum for $\int_{\RR^d} \Abs{\nabla\ft\phi}$,
which is finite, hence $\lim_{t\to\infty} S_2 = 0$.

Since $\epsilon$ is arbitrary the proof is complete.
\Qed

\begin{remark}\label{rem:point-support}
The same proof as that of Theorem \ref{th:at-zero} shows that,
if
$$
\mu = \sum_{\lambda\in\Lambda} c_\lambda \delta_\lambda,
$$
with $\Abs{c_\lambda} \le C$, $\Lambda$ is of density $0$ and
the tempered distribution $\ft\mu$ is locally a measure
in the neighborhood of some point $a\in\RR^2$, then we 
have $\ft\mu(\Set{a}) = 0$.
\end{remark}

\begin{theorem}\label{th:point-implies-measure}
Suppose that the multiset $\Lambda\subset\RR^d$ has density
uniformly bounded by $\rho$
and that, for some point $a\in\RR^d$ and $R>0$,
$$
\supp{\ft{\delta_\Lambda}} \cap B_R(a) = \Set{a}.
$$
Then, in $B_R(a)$, we have $\ft{\delta_\Lambda} = w \delta_a$, for
some $w \in \CC$ with $\Abs{w} \le \rho$.
\end{theorem}  

\noindent{\bf Proof.}
It is well known that the only tempered distributions supported
at a point $a$ are finite linear combinations of the derivatives
of $\delta_a$.
So we may assume that, for $\phi\in C^\infty(B_R(a))$,
\beql{point-distr}
\ft{\delta_\Lambda}(\phi) =
 \sum_\alpha c_\alpha (D^\alpha \delta_a) (\phi) =
 \sum_\alpha (-1)^\Abs{\alpha} c_\alpha D^\alpha \phi(a),
\eeq
where the sum extends over all values of the multiindex
$\alpha=(\alpha_1,\ldots,\alpha_d)$ with
$\Abs{\alpha} = \alpha_1+\cdots+\alpha_d \le m$ (the finite degree)
and $D^\alpha = \partial_1^{\alpha_1}\cdots\partial_d^{\alpha_d}$
as usual.

We want to show that $m=0$. Assume the contrary and let $\alpha_0$
be a multiindex that appears in \eqref{point-distr}
with a non-zero coefficient and has $\Abs{\alpha_0} = m$.
Pick a smooth function $\phi$
supported in a neighborhood of $0$ which is such that for each
multiindex $\alpha$ with $\Abs{\alpha} \le m$ we have
$D^\alpha \phi (0) = 0$ if $\alpha \neq \alpha_0$ and
$D^{\alpha_0} \phi(0) = 1$. (To construct such a $\phi$,
multiply the polynomial $(1/\alpha_0!) x^{\alpha_0}$ with a
smooth function supported in a neighborhood of $0$, which
is identically equal to $1$ in a neighborhood of $0$.)

For $t\to\infty$ let $\phi_t(x) = \phi(t(x-a))$.
Equation \eqref{point-distr} then gives that
\beql{one-side}
\ft{\delta_\Lambda}(\phi_t) = t^m (-1)^m c_{\alpha_0}.
\eeq
On the other hand, using
$$
\FT{\phi(t(x-a))}(\xi) = e^{-2\pi i\inner{a}{\xi/t}} t^{-d}
\ft\phi(\xi/t),
$$
we get
\beql{other-side}
\ft{\delta_\Lambda}(\phi_t) =
 \sum_{\lambda\in\Lambda} e^{-2\pi i\inner{a}{\lambda/t}} t^{-d}
\ft\phi(\lambda/t).
\eeq
Notice that \eqref{other-side} is a bounded quantity as $t\to\infty$
by a proof similar to that of Theorem \ref{th:at-zero},
while \eqref{one-side} increases like $t^m$, a contradiction.

Hence $\ft{\delta_\Lambda} = w \delta_a$ in a neighborhood of $a$.
The proof of Theorem \ref{th:at-zero} again gives that $\Abs{w} \le
\rho$.
\Qed

We are now ready to prove the result \cite{kolountzakis:polygons} that finishes the argument.
\begin{theorem}\label{th:periodic-set}
{\bf (Kolountzakis, 2000)}\\
Suppose that $\Lambda \subset \RR^2$ is a discrete
multiset of uniformly bounded density and that
$$
\ft{\delta_\Lambda} = \FT{\sum_{\lambda\in\Lambda} \delta_\lambda}
$$
is locally a measure with
$$
\Abs{\ft{\delta_\Lambda}}(B_R(0)) \le C R^2,
$$
for some positive constant $C$ and $R \ge 1$.
Assume also that $\ft{\delta_\Lambda}$ has discrete support.
Then $\Lambda$ is a finite union of translated lattices.
\end{theorem}

\noindent{\bf Proof.}
Define the sets (not multisets)
$$
\Lambda_k = \Set{\lambda \in \Lambda:\ \mbox{$\lambda$
 has multiplicity $k$}}.
$$
By Meyer's Theorem \ref{th:meyer} (applied for the base set of the multiset $\Lambda$
with the coefficients $c_\lambda$ equal to the corresponding
multiplicities) each of the $\Lambda_k$ is in the coset ring
of $\RR^2$.

By Theorem \ref{th:coset-ring} it follows that the discrete set $\Lambda_k$
can be constructed from lattices in $\RR^2$ (two-dimensional, one-dimensional or points)
using finitely many operations and one shows easily that the set $\Lambda_k$ has the form
\beql{form}
\Lambda_k ~ = ~ \left( \bigcup_{j=1}^J A_j ~ \setminus ~ (B_1^{(j)} \cup 
\cdots \cup
 B_{n_j}^{(j)})\right)~ \cup ~ \bigcup_{l=1}^L L_l ~ \bigtriangleup
~ F,
\eeq
where $A_1,\ldots,A_J$ are $2$-dimensional translated lattices,
$L_l$ and
$B_i^{(j)}$ are
$1$-dimensional translated lattices and $F$ is a finite set ($J, L
\ge 0$).
The lattices $A_j$ may be assumed to be have pairwise intersections of
dimension at most $1$.

We may thus write
\beql{identity}
\Lambda_k = A \bigtriangleup B,
\eeq
with $A = \bigcup_{j=1}^J A_j$, where the $2$-dimensional
translated lattices $A_j$ have pairwise intersections of
dimension at most $1$, and $\dens B = 0$.

Hence
$$
\delta_{\Lambda_k} = \sum_{j=1}^J \delta_{A_j} + \mu,
$$
where $\mu = \sum_{f \in F} c_f \delta_f$, $\dens F = 0$ and
$\Abs{c_f} \le C(J)$.
The set $F$ consists of $B$ and all points contained in at least
two of the $A_j$.

Combining for all $k$, and reusing the symbols
$A_j$, $\mu$ and $F$, we get
$$
\delta_{\Lambda} = \sum_{j=1}^J \delta_{A_j} + \mu.
$$
But $\ft{\delta_\Lambda}$ and $\sum_{j=1}^J \ft{\delta_{A_j}}$ are
both (by the assumption and the Poisson Summation Formula) discrete
measures, and so is therefore $\ft\mu$.
However $\dens F = 0$ and the boundedness of the coefficients $c_f$
implies that $\ft\mu$ has no point masses
(see Remark \ref{rem:point-support}),
which means that $\ft\mu = 0$ and so is $\mu$.
Hence $\delta_{\Lambda} = \sum_{j=1}^J \delta_{A_j}$, or
$$
\Lambda = \bigcup_{j=1}^J A_j,\ \ \ \mbox{as multisets}.
$$
\Qed

Last, observe that the support of $\ft{\delta_\Lambda}$ is contained
in the intersection of two grids of the type shown in Theorem \ref{th:polygon-zero-set},
and has therefore (remember it's a discrete set) bounded density.
This proves that
$\Abs{\ft{\delta_\Lambda}}(B_R(0)) \le C R^2$
and we can invoke Theorem \ref{th:periodic-set}.

\pagebreak
\section{Lecture 2: Problems of lattice tiling.}

Here we will examine several lattice tiling problems. The study of lattice tilings
in Fourier space is particularly simple as explained in \S \ref{sec:lattice}
\begin{quotation}
$f+\Lambda$ is a tiling if and only if $\ft f$
vanishes on the dual lattice $\Lambda^*$,
except at zero.
\end{quotation}
The study of lattice tiling does not involve at all distributions
which are not measures.
The Fourier Analysis involved is nothing more than the usual multi-dimensional
Fourier Series plus a change of variable to go from the integer to the
arbitrary lattice.

\subsection{A new equivalent form of a theorem of Haj\'os}
\label{sec:new-form-of-hajos}

Let us start by quoting a well known theorem of Minkowski in the Geometry
of Numbers.
\begin{theorem}\label{th:minkowski-linear-forms}
{\bf (Minkowski, ca. 1900)}\\
Let $A \in GL(d, \RR)$ have $\det A = 1$.
Then there is $x\in\ZZ^d\setminus\Set{0}$ with $\Linf{Ax} \le 1$.
\end{theorem}

\noindent{\bf Proof.}
Let $\Lambda = A\ZZ^d$ and $U = [-{1\over 2}, {1\over 2}]$. We want to show
that $\Lambda \cap (2U)$ contains something besides $0$.
Suppose, on the contrary, that $\Lambda \cap (2U) = \Set{0}$.
Then, there is $\epsilon>0$ such that for
$$
U_\epsilon = [-{1\over 2}-\epsilon, {1\over 2}+\epsilon]
$$
we have $\Lambda \cap (2U_\epsilon) = \Set{0}$.
We can rewrite this as
$$
(\Lambda - \Lambda) \cap (U_\epsilon - U_\epsilon) = \Set{0},
$$
which means that the copies $U_\epsilon+\lambda$, $\lambda \in \Lambda$,
are disjoint (we have a {\em packing}).
But $\dens \Lambda = 1$ and $\Abs{U_\epsilon} > 1$, 
which is a contradiction, according to Lemma \ref{lm:packing-implies-density-bound}.
\Qed

The following theorem of Haj\'os \cite{hajos} proved a conjecture of Minkowski some
forty years after it was posed.
This conjecture concerned the case when one could have a strict inequality in Theorem
\ref{th:minkowski-linear-forms}.
\begin{theorem}\label{th:hajos-forms}
{\bf (Haj\'os, 1941)}\\
Let $A \in GL(d, \RR)$ have $\det A = 1$.
Then there is $x\in\ZZ^d$ with $\Linf{Ax} < 1$
unless $A$ has an integral row.
\end{theorem}

Haj\'os actually worked on the following equivalent form
of the Minkowski conjecture, which involves lattice tilings by a cube.
This form was already known to Minkowski and most results on Minkowski's
conjecture leading up to Haj\'os's eventual proof have used
this form.
\begin{theorem}\label{th:hajos-tilings}
If $Q = [-1/2, 1/2]^d$ is a cube of unit volume in $\RR^d$,
$\Lambda\subset\RR^d$ is a lattice, and
$$
\RR^d = Q + \Lambda
$$
is a lattice tiling of $\RR^d$ then there are two cubes in the tiling
that share a $(d-1)$-dimensional face.
In other words,
for some $i=1,\ldots,d$, the standard basis vector
$e_i = (0,\ldots,0,1,0,\ldots,0)^\top \in \Lambda$.
\end{theorem}

Keller \cite{keller} conjectured that the same is true even without the
lattice assumption. That is, Keller conjectured that in any tiling of Euclidean
space by translates of a cube there are two cubes in the tiling which share
a $(d-1)$-dimensional face.
This is indeed true up to dimension $6$ but was disproved by Lagarias and Shor
\cite{lagarias-shor} for $d\ge 10$.
Mackey \cite{mackey} recently proved that Keller's conjecture
is false in dimension 8 and higher.
The only remaining case for which Keller's conjecture is still open is $d=7$.

\noindent{\bf Theorem \ref{th:hajos-forms} $\Longrightarrow$ Theorem
\ref{th:hajos-tilings}.}

\noindent
Let $\Lambda = A \ZZ^d$ with $\det A =1$, $Q+\Lambda = \RR^d$.
Then, either there is a non-zero $\Lambda$-point in the interior
of $2Q$ or $A$ has an integral row.
The first cannot happen because of the tiling assumption.
Therefore $a_{ij} \in \ZZ$ for some $i$ and for all $j$.
Again because of tiling it follows that $\gcd(a_{i1},\ldots,a_{i,d})=1$.
Otherwise the $i$-th coordinates of all $\Lambda$-points would be multiples
of $G=\gcd(a_{i1},\ldots,a_{i,d})>1$, which is impossible (there would be gaps in
the tiling).
Let $\RR^{d-1}$ be the subspace spanned by all $e_j$, $j\neq i$, and
define $\Lambda' = \Lambda \cap \RR^{d-1}$ and
$Q' = Q \cap \RR^{d-1}$.
It follows that $\RR^{d-1} = \Lambda' + Q'$ is a tiling
of $\RR^{d-1}$. By induction then $\Lambda'$ contains some
vector of the standard basis and so does $\Lambda$.
\Qed

\noindent{\bf Theorem \ref{th:hajos-tilings} $\Longrightarrow$ Theorem
\ref{th:hajos-forms}.}

Theorem \ref{th:hajos-tilings} easily implies the
seemingly stronger statement that, if $A\ZZ^d + Q = \RR^d$ is a
tiling then, after a permutation of the coordinate axes, the matrix $A$
takes the form
\beql{matrix-form}
\left(
\begin{array}{ccccc}
1 & 0 & 0 & \ldots & 0 \\
a_{2,1} & 1 & 0 & \ & 0 \\
\ldots & \ldots & \ldots & \ldots & \ldots\\
a_{d,1} & \ldots &\ldots &\ldots & 1
\end{array}
\right)
\eeq
Using this remark, if $A\ZZ^d \cap (-1,1)^d = \Set{0}$ we get, since
$\det A = 1$, that $A\ZZ^d + Q = \RR^d$ and, therefore, $A$ is
(after permutation of the coordinate axes) of the type
\eqref{matrix-form}, and thus has an integral row (and this property
is preserved under permutation similarity).
\Qed

We now prove that the following is equivalent
to Theorems \ref{th:hajos-forms} and \ref{th:hajos-tilings} \cite{kolountzakis:notched}.

\begin{theorem}\label{th:3rd}
{\bf (Kolountzakis, 1998)}\\
Let $B \in GL(d, \RR)$ have $\det B = 1$ and the property that
for all $x\in\ZZ^d\setminus\Set{0}$ some coordinate
of the vector $Bx$ is a non-zero integer.
Then $B$ has an integral row.
\end{theorem}

\begin{open-problem}
Prove this combinatorial statement directly, thereby obtaining
a new proof of the Minkowski Conjecture.
\end{open-problem}

\begin{remark}
One might think that Theorem \ref{th:3rd} can be proved
equivalent directly to Theorem \ref{th:hajos-forms}, which
it resembles most.
It is, indeed, clear that Theorem \ref{th:hajos-forms} implies
Theorem \ref{th:3rd}.
However, the proof that is given here is that of the equivalence
of Theorems \ref{th:3rd} and \ref{th:hajos-tilings}.
I do not know of a more direct proof of the fact that
Theorem \ref{th:3rd} implies Theorem \ref{th:hajos-forms}.
\end{remark}

We shall need the following simple lemma.
\begin{lemma}\label{lm:matrix}
Let $A \in GL(d, \RR)$ be a non-singular matrix.
The lattice $A^{-\top}\ZZ^d$ contains the basis vector $e_i$ 
if and only if the $i$-th row of $A$ is integral.
\end{lemma}

\noindent{\bf Proof.}
Without loss of generality assume $i=1$.

If $e_1 \in A^{-\top}\ZZ^d$ then $e_1 = A^{-\top} x$ for some
$x\in\ZZ^d$.
Therefore, for all $y\in\ZZ^d$ we have
$$
(Ay)_1 = e_1^\top A y = x^\top A^{-1} A y = x^\top y \in \ZZ.
$$
It follows that $(Ay)_1 \in \ZZ$ for all $y \in \ZZ^d$ and the
first row of $A$ is integral.

Conversely, if the first row of $A$ is integral, then,
for all $y\in\ZZ^d$
$$
\ZZ \ni (Ay)_1 = x^\top y,
$$
where $A^{-\top}x = e_1$ ($x\in\RR^d$).
It follows that $x\in\ZZ^d$ and $e_1 \in A^{-\top}\ZZ^d$.
\Qed

\noindent{\bf 
Proof of the equivalence of Theorems \ref{th:hajos-tilings}
and \ref{th:3rd}.
}

Let $f(x) = \One{x \in Q}$ be the indicator function
of the unit-volume cube $Q = [-1/2, 1/2]^d$.
A simple calculation shows that
\beql{FT-of-f}
\widehat{f}(\xi)  = \prod_{j=1}^d {\sin \pi \xi_j \over \pi \xi_j},
\eeq
so that
\beql{Z-of-f}
{\cal Z} := \Set{\widehat{f} = 0} = \Set{\xi\in\RR^d:\ \mbox{some $\xi_j$ is
	a non-zero integer}}.
\eeq
Therefore, if $\Lambda = B^{-\top}\ZZ^d$ then (since $\Lambda$ has
volume $1$)
$$
Q+\Lambda = \RR^d \Longleftrightarrow
	\Lambda^*\setminus\Set{0} \subseteq {\cal Z},
$$
where $\Lambda^* = B\ZZ^d$.
In words, $Q$ tiles with $\Lambda=B^{-\top}\ZZ^d$ if and only if
for every $x\in\ZZ^d\setminus\Set{0}$ the vector $Bx$ has
some non-zero integral coordinate.

\noindent{\bf Theorem \ref{th:hajos-tilings} $\Longrightarrow$
 Theorem \ref{th:3rd}}.

\noindent
Suppose $x\in\ZZ^d\setminus\Set{0}$ implies some $(Bx)_i\in\ZZ\setminus\Set{0}$.
Then $Q+\Lambda = \RR^d$ and from Theorem \ref{th:hajos-tilings},
say, $e_1 \in \Lambda$, which, from Lemma \ref{lm:matrix},
implies that the first row of $B$ is integral.

\noindent{\bf Theorem \ref{th:3rd} $\Longrightarrow$
  Theorem \ref{th:hajos-tilings}}.

\noindent
Assume $Q+\Lambda = \RR^d$. It follows that for every
$x\in\ZZ^d\setminus\Set{0}$ the vector $Bx$ has some non-zero integral
coordinate. By Theorem \ref{th:3rd} $B$ must have an integral row,
which, by Lemma \ref{lm:matrix}, implies that some $e_i \in \Lambda$.
\Qed

\subsection{Tilings by notched and extended cubes}
\label{sec:notched-extended}

\begin{figure}[h] \centering
\input{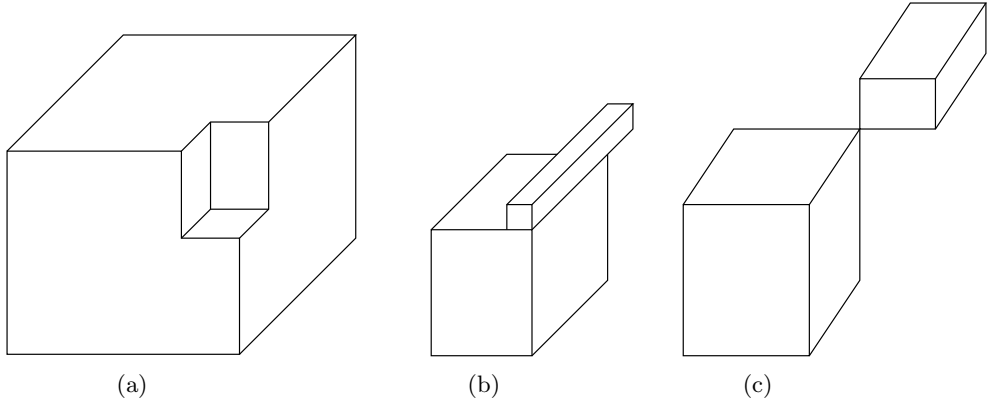}
\caption{These shapes admit lattice tilings}
\ \\ \hrule
\label{fig:cubes}
\end{figure}

\noindent
In this section we prove that some simple shapes (like those in Figure
\ref{fig:cubes}) admit lattice tilings.

\subsubsection{The notched cube}
We consider first the unit cube
$$
Q = \left[-{1\over2},{1\over2}\right]^d
$$
from whose corner (say in the positive orthant) a rectangle
$R$ has been removed with sides-lengths $\delta_1,\ldots,\delta_d$
($0\le\delta_j\le1$).
That is, we consider the ``notched cube'':
$$
N = Q \setminus R
$$
where
$$
R = \prod_{j=1}^d \left[{1\over2}-\delta_j,{1\over2}\right].
$$
It is shown in Figure \ref{fig:cubes} (a).

We give a new \cite{kolountzakis:notched},
Fourier-analytic, proof of the following result of Stein \cite{stein}.
\begin{theorem}\label{th:notched}
{\bf (Stein, 1990)}\\
The notched cube $N$ admits a lattice tiling of $\RR^d$.
\end{theorem}

After a simple calculation we obtain
\beql{cube-FT}
\widehat{\chi_N}(\xi) = \prod_{j=1}^d {\sin\pi\xi_j \over \pi\xi_j} -
  F(\xi) \prod_{j=1}^d {\sin\pi\delta_j\xi_j \over \pi\xi_j},
\eeq
where $F(\xi) = \exp(\pi i K(\xi))$ with
\beql{K}
K(\xi) = \sum_{j=1}^d (\delta_j-1)\xi_j.
\eeq
Using \eqref{lattice-tiling-in-fourier-space} it is enough
to exhibit a lattice $\Lambda \subset \RR^d$, of volume
equal to
$$
\Abs{N} = 1-\delta_1\cdots\delta_d,
$$
such that
$\widehat{\chi_N}$ vanishes on $\Lambda^*\setminus\Set{0}$.

\subsubsection{Lattices in the zero-set}
\label{sec:zeros}
We define the lattice $\Lambda^*$ as those points $\xi$ for
which
\begin{eqnarray}\label{system}
\xi_1 - \delta_2\xi_2 &=& n_1,\nonumber\\
\xi_2 - \delta_3\xi_3 &=& n_2,\nonumber\\
 \ldots & & \\
\xi_d - \delta_1\xi_1 &=& n_d,\nonumber
\end{eqnarray}
for some $n_1,\ldots,n_d \in \ZZ$.
That is, $\Lambda^* = A^{-1} \ZZ^d$, where
\beql{matrix}
A = \left(
\begin{array}{cccccc}
1 & -\delta_2 & & & & \\
 & 1 & -\delta_3 & & & \\
 & & & \ddots & & \\
 & & & & 1 & -\delta_d \\
-\delta_1 & & & & & 1
\end{array} \right).
\eeq
Therefore $\Lambda = A^\top\ZZ^d$ and the volume of $\Lambda$
is equal to $\Abs{\det A}$.
Expanding $A$ along the first column we get easily
that $\det A = 1-\delta_1\cdots\delta_d$,
which is the required volume.

We now verify that $\widehat{\chi_N}$ vanishes on
$\Lambda^*\setminus\Set{0}$.

Assume that $0 \neq \xi \in \Lambda^*$.
Adding up the equations in \eqref{system}
we get
$$
K = K(\xi) = -(n_1+\cdots+n_d).
$$
If all the coordinates of $\xi$ are non-zero
we can write
\beql{tmp1}
\widehat{\chi_N}(\xi) = {1\over\pi^d\xi_1\cdots\xi_d} \left(
 \prod_{j=1}^d\sin\pi\xi_j - (-1)^{K(\xi)} \prod_{j=1}^d\sin\pi\delta_j\xi_j \right).
\eeq
Observe from \eqref{system} that
$$
\sin\pi\xi_j = (-1)^{n_j} \sin\pi\delta_{j+1}\xi_{j+1},
$$
where the subscript arithmetic is done modulo $d$,
 from which we get $\widehat{\chi_N}(\xi) = 0$,
since the factors in the two terms of \eqref{tmp1}
match one by one.

It remains to show that $\widehat{\chi_N}(\xi) = 0$ even when
$\xi$ has some coordinate equal to $0$, say $\xi_1=0$.

Consider the numbers $\xi_1,\ldots,\xi_d$ arranged in a cycle
and let
$$
I = \Set{\xi_m, \xi_{m+1},\ldots,\xi_1,\ldots,\xi_{k-1},\xi_k}
$$
be an interval around $\xi_1$ which is maximal with the property
that all its elements are $0$.
Then $\xi_{m-1}\neq 0$ and $\xi_{k+1}\neq 0$
and from \eqref{system} we get
\beql{eqs}
\xi_{m-1} - \delta_m\xi_m = n_m
\ \ \mbox{and}\ \ 
\xi_k - \delta_{k+1}\xi_{k+1} = n_k.
\eeq
We deduce that $n_m$ and $n_k$ are both non-zero
and therefore that
$\xi_{m-1}$ and $\delta_{k+1}\xi_{k+1}$ are both
non-zero integers
and $\sin \pi\xi_{m-1} = \sin \pi\delta_{k+1}\xi_{k+1} = 0$.
This means that both terms in \eqref{cube-FT} vanish and
so does $\widehat{\chi_N}(\xi)$.
\Qed

So we proved that for the lattice $\Lambda = A^\top\ZZ^d$,
where $A$ is defined in \eqref{matrix},
we have $N+\Lambda = \RR^d$.
Clearly, if $\sigma$ is a cyclic permutation
of $\Set{1,\ldots,d}$ and
if instead of the matrix $A$ we have the matrix $A'$ whose $i$-th
row has $1$ on the diagonal, $-\delta_{\sigma i}$ at column $\sigma i$ and
$0$ elsewhere, we get again a lattice tiling with
the lattice $(A')^\top\ZZ^d$.
Stein \cite{stein} as well as Schmerl \cite{schmerl}
have shown that these $(d-1)!$ lattice tilings of the notched
cube (one for each cyclic permutation of $\Set{1.,\ldots,d}$)
are all non-isometric when the side-lengths $\delta_j$ are
all distinct.

A deeper result of Schmerl \cite{schmerl}
is that there are no other translational tilings of the notched cube,
lattice or not.
This is something that cannot apparently be proved with the Fourier Analysis
approach.

\subsubsection{Extended cubes}
\label{sec:ext}
Let us now allow the parameters $\delta_1,\ldots,\delta_d$
to take on any non-zero real value
subject only to the restriction
\beql{not-one}
\delta_1\cdots\delta_d \neq 1,
\eeq
and let the function $\varphi(\xi)$
be equal to the right-hand side of \eqref{cube-FT}.
Let again the matrix $A$ be defined by \eqref{matrix}
and $\Lambda = A^\top\ZZ^d$ as before.
We have again $\det A = 1-\delta_1\cdots\delta_d$.

The calculations we did in \S\ref{sec:zeros} show
that $\varphi$ vanishes on $\Lambda^*\setminus\Set{0}$, hence,
if $\check{\varphi}$ is the inverse Fourier Transform of $\varphi$, $\check{\varphi}$
tiles $\RR^d$ with $\Lambda$ and weight
\beql{weight}
{\varphi(0) \over \Abs{1-\delta_1\cdots\delta_d}} =
	{\rm sgn}(1-\delta_1\cdots\delta_d),
\eeq
where ${\rm sgn}(x)=\pm1$ is the sign of $x$.

The function $\check{\varphi}$ is given by
\beql{IFT}
\check{\varphi}(x) = \chi_Q(x) - {\rm sgn}(\delta_1\cdots\delta_d)\psi(x),
\eeq
where
\beql{psi}
\psi(x) = 
 \chi_Q\left({x_1-(1-\delta_1)/2 \over \Abs{\delta_1}},\ldots,
	     {x_d-(1-\delta_d)/2 \over \Abs{\delta_d}} \right).
\eeq
Notice that $\psi(x)$ is the indicator function
of a rectangle $R = R(\delta_1,\ldots,\delta_d)$ with
side-lengths $\Abs{\delta_1},\ldots,\Abs{\delta_d}$
centered at the point
\beql{ptP}
P = \left({1\over2},\ldots,{1\over2}\right) -
	{1\over2}\left(\delta_1,\ldots,\delta_d\right).
\eeq
The rectangle $R$ intersects the interior of $Q$ only in the case
$\delta_1>0,\ldots,\delta_d>0$ and when this happens
$\check{\varphi}$ is an indicator function only if we
also have $\delta_1\le1,\ldots,\delta_d\le1$, which is the case
of the notched cube that we examined in \S\ref{sec:zeros}.

Otherwise (not all the $\delta$s are non-negative)
$\check{\varphi}$ is an indicator function only when
${\rm sgn}(\delta_1\cdots\delta_d) = -1$, i.e.,
the number of negative $\delta$s is odd.
In this case we have that
$$
\check{\varphi} = \chi_{Q\cup R}
$$
and from \eqref{weight} we get that $Q\cup R$
tiles with $\Lambda$ and weight $1$.
We can now prove the following \cite{kolountzakis:notched}.
\begin{theorem}\label{th:ext}
{\bf (Kolountzakis, 1998)}\\
Let $Q$ and $R$ be two axis-aligned rectangles in $\RR^d$
with sides of arbitrary length and disjoint interiors.
Assume also that $Q$ and $R$ have a vertex $K$ in common
and intersection of odd codimension.

Then $Q\cup R$ admits a lattice tiling of $\RR^d$ of weight $1$.
\end{theorem}
For example, the extended cubes shown in Figure \ref{fig:cubes} (b),(c)
admit lattice tilings of $\RR^3$, as the corresponding codimensions are
$1$ and $3$.

\noindent{\bf Proof.}
After a linear transformation we can assume
that $Q = [-1/2,1/2]^d$, that $Q$ and $R$ share the vertex
$K = (1/2,\ldots,1/2)$ and that
$Q\cap R$ has codimension $k$ (an odd number) and
$$
Q \cap R \ \subseteq\ \Set{x\in Q:\ x_1=\cdots=x_k={1\over2}}.
$$

Let the side-lengths of $R$ be $\gamma_1,\ldots,\gamma_d>0$.
Define
$$
\delta_j = \left\{\begin{array}{ll}
 -\gamma_j, & \mbox{if $1\le j \le k$},\\
 \gamma_j, & \mbox{if $k+1\le j \le d$.}
\end{array}\right.
$$
It follows that, with this assignment for the $\delta_j$,
the indicator function of $R$ is equal to the function
$-{\rm sgn}(\delta_1\cdots\delta_d)\psi(x)$
of \eqref{IFT} and tiling follows from the previous discussion.
\Qed

Most likely the extended cubes with an intersection of even codimension
do not tile, at least not for general side-lengths.
This is clear in dimension two and it is conceivable that some
combinatorial argument could easily show this in any dimension.
The Fourier Analysis approach does not seem to be very helpful
when one tries to disprove that something is a translational tile.

\begin{open-problem}
In the setting of Theorem \ref{th:ext} prove that if the
codimension is even then the set $Q \cup R$ is not a tile.
\end{open-problem}

\subsection{The Steinhaus tiling problem}

\subsubsection{The original, two-dimensional case}

\noindent
Steinhaus \cite[problem 59]{moser} asked whether there is
a planar set $S$ which, no matter how translated and rotated,
always contains exactly one point with integer coordinates.
\begin{definition}\label{def:steinhaus}
{\bf (Steinhaus property)}\\
A set $S\subset\RR^2$ has the {\em Steinhaus property}
if for every $x\in\RR^2$ and for every rotation
$$
A_\theta = \left(\begin{array}{ll} \cos\theta & -\sin\theta\\
                             \sin\theta & \cos\theta \end{array}\right)
$$
we have
\begin{equation}\label{eq:steinhaus1}
\#\left( \ZZ^2 \cap (A_\theta S + x) \right) = 1,
\end{equation}
where $A_\theta S + x = \{A_\theta s + x :\ s \in S\}$.
\end{definition}
Sierpi\'nski \cite{sierpinski} first proved that a set which
is bounded and either open or closed cannot have the Steinhaus property.
Croft \cite{croft} and Beck \cite{beck} proved the same
of any set which is bounded and measurable.
(Croft's approach is more direct and geometric. Beck is using
Fourier Analysis.)
Ciucu \cite{ciucu} shows that any Steinhaus set must have empty
interior, without assuming boundedness.
Several variations of the problem have been investigated
by Komj\'ath \cite{komjath} from a rather different point a view, where one
places a different subgroup of the plane in place of $\ZZ^2$.

Very recently it was shown by Jackson and Mauldin \cite{jackson-mauldin}
that Steinhaus sets do indeed exist. But the construction there does not
furnish measurable such sets and it is precisely under the assumption
of measurability that we study the existence problem for Steinhaus sets here,
using Fourier Analysis.

To begin, notice that the question of Steinhaus can be rephrased as follows:
\begin{itemize}
\item[(a)] Is there a set $E$ which tiles the plane
if translated at any rotated copy of $\ZZ^2$?
\item[(b)] Or, is there a common set of coset representatives (fundamental domain)
of all groups $R_\theta\ZZ^2$
in the group $\RR^2$?
\end{itemize}
We only care for measurable Steinhaus sets (if they exist) so tiling, above,
is to be interpreted in the almost everywhere sense, as it is normally interpreted
throughout this survey.

As first noticed by Beck \cite{beck}, the Steinhaus question in the form (a), above,
is equivalent to asking if there exists
a measurable set $E \subseteq \RR^2$, of measure $1$,
such that the Fourier Transform of its indicator function vanishes on all circles of
the plane which are centered at the origin and pass through some point of the
integer lattice $\ZZ^2$.
This is so since for a set to have the Steinhaus property it must tile the plane
when translated by any rotation of $\ZZ^2$ (this alone implies
of course that $\Abs{E}=1$). These sets are lattices, hence this
is equivalent to $\ft{\chi_E}$ vanishing on all these lattices, which
are self-dual. 
The union of these rotated lattices
is precisely the set of circles mentioned above.
We state this as a Theorem.
\begin{theorem}\label{th:steinhaus-equivalence}
A measurable set $E \subseteq \RR^2$ is simultaneously a tile for
all rotations of $\ZZ^2$ if and only if it has measure $1$ and
its Fourier Transform $\ft{\chi_E}$ vanishes on all circles with
center at the origin and radius of the form $\sqrt{m^2+n^2}$,
with $m,n \in \NN$, not both $0$.
\end{theorem}
It is now easy to see
that such sets cannot be bounded, if they exist.
Indeed, the restriction onto any line $L$ through $0$ of $\ft{\chi_E}$
is nothing but the one-dimensional Fourier Transform of the function
$\chi_E$ projected onto $L$, i.e., of the function
$$
f(t) = \int_{L^\perp} \chi_E(tu + s)~ds,
$$
where $u$ is a unit vector on $L$ and $L^\perp$ is the line through $0$ which is
orthogonal to $L$.
But if $E$ is bounded the function $f(t)$ has compact support,
hence $\ft{\chi_E}(t u)$ is an entire function of exponential type,
and, as such, it should have at most $C\cdot R$ zeros in the
interval $(-R, R)$, where $C>0$ is a constant.
(See the discussion in \S\ref{sec:structure-in-dim-1}.)
However, the number of zeros of $\ft{\chi_E}(t u)$ is twice the number
of circles out to radius $R$, or, in other words, twice the number of
integers expressible as a sum of two integer squares and of size up to $R^2$.
But this number is almost quadratic in $R$.
It is a well known result of Landau \cite{fricker} that it is
$\sim c R^2 \log^{-1/2}R$.

With a more careful and quantitative approach along similar lines,
but not using entire functions, it
was then proved by the author \cite{kolountzakis:steinhaus}
that any set $E$ with the Steinhaus property must be large at infinity:
$$
\int_E \Abs{x}^\alpha ~dx = \infty,\ \ \ \mbox{for any $\alpha>{10 \over 3}$}.
$$

With much more care it was obtained in \cite{kolountzakis-wolff}
by the author and Tom Wolff that
\begin{theorem}\label{th:kw-two-dim}
{\bf (Kolountzakis and Wolff, 1997)}\\
If $E\subseteq\RR^2$ is a measurable Steinhaus set then $\int_E\Abs{x}^\alpha=\infty$,
for all $\alpha>46/27$.
\end{theorem}
The number $46/27$ comes from the best known estimate known for the {\em circle problem}.
This is the problem where one asks for the best upper estimates in the error
term $E(R)$ (as $R\to\infty$) in the expression
$$
N(R) = \pi R^2 + E(R),
$$
where $N(R)$ is the number of integer lattice points
in the disk $\Set{\Abs{x}\le R}\subseteq\RR^2$.
Even if the conjectured best possible upper bound $E(R) = O(R^{1/2 + \epsilon})$ gets proved
the estimate for the Steinhaus tiling problem in Theorem \ref{th:kw-two-dim} would
only become true for all $\alpha > 1$. So it appears that if one is going to disprove
the existence of measurable Steinhaus sets in dimension two one needs some rather different approach.

This seems to be the state of knowledge for the two-dimensional case.

\subsubsection{The problem in dimension $d\ge3$}

\noindent
The Steinhaus problem generalizes very naturally to any dimension. One
asks for a set $E \subseteq \RR^d$ such that no matter what orthogonal
linear transformation you apply to it, it still tiles $\RR^d$ when translated
by $\ZZ^d$. With precisely the same reasoning as before, one is looking for a
measurable set of measure $1$ such that the Fourier Transform of its indicator function
vanishes on all spheres centered at the origin that contain some integer
lattice point.

It is because of the fact that we know precisely which numbers are representable
as sums of three squares that the following result \cite{kolountzakis-wolff} holds.
\begin{theorem}\label{th:steinhaus-functions-are-continuous}
{\bf (Kolountzakis and Wolff, 1997)}\\
If $f\in L^1(\RR^d)$, $d\ge 3$, and $\ft f$ vanishes on all
spheres centered at the origin through some lattice point,
then $f$ is a.e.\ equal to a continuous function.

In particular, there are no measurable Steinhaus sets in dimension
$d\ge 3$.
\end{theorem}

Here we show an alternative way \cite{kolountzakis-papadimitrakis:quad} 
of proving that there are no sets
with the Steinhaus property in dimension $d\ge 3$.
We emphasize though that Theorem \ref{th:steinhaus-functions-are-continuous}
is much stronger than Theorem \ref{th:quad} given below.
See also some related results of Mauldin and Yingst \cite{mauldin-yingst}.
\begin{theorem}\label{th:quad}
{\bf (Kolountzakis and Papadimitrakis, 2000)}\\
There are no measurable Steinhaus sets in dimension $d\ge 3$.
\end{theorem}

\noindent{\bf Proof.}
In any dimension $d$ write ${\cal B}$ for the union of all spheres
centered at the origin that go through at least one lattice point.
The point $0$ is included in ${\cal B}$.

Assume from now on that the set $E$ is a Steinhaus set in dimension $d$.

Suppose now that we can find a lattice $\Lambda^* \subset {\cal B}$
with $\det\Lambda^*$ not an integer.
Since $\ft{\chi_E}$ vanishes on $\Lambda^*\setminus\Set{0}$ it
follows that $E + \Lambda$ is a tiling at level
$\ell = \Abs{E}\times\dens\Lambda = 1 \times \det\Lambda^*$,
which is not an integer.
This is a contradiction as, obviously, any set may only
tile at an integral level.

Looking at the quadratic form $\inner{A^\top Ax}{x}$
for each lattice $\Lambda^* = A\ZZ^d$
we summarize the above observations in the following lemma
\begin{lemma}\label{lm:quad}
If there exists a positive definite quadratic form
$Q(x) = Q(x_1,\ldots,x_d) = \inner{Bx}{x}$ such that for all integral
$x_1,\ldots,x_d$ its value is the sum of $d$ integer squares,
and the determinant of $Q$, $\det B$, is not the square of an integer,
then there are no Steinhaus sets in dimension $d$.
\end{lemma}

\noindent
{\bf The case $d\ge 4$:}\\
Consider the symmetric $4\times 4$ matrix $B$ with $1$ on the diagonal
and $1/2$ everywhere else.
The matrix $B$ is positive definite (its eigenvalues are
$1/2$, $1/2$, $1/2$ and $5/2$) and its determinant is
$5/16$.
It defines the quadratic form
$$
Q(x) = Q(x_1, \ldots, x_4) = \inner{Bx}{x} = \sum_{i=1}^4 x_i^2 +
        \sum_{i>j} x_i x_j,
$$
which is obviously integer valued and has non-square determinant.
Furthermore, every non-negative integer may be written as a sum
of four squares (Lagrange).
From Lemma \ref{lm:quad} it follows that there are no
Steinhaus sets for $d=4$. We easily see that this extends to all higher dimensions
by taking as our matrix the identity in one corner of which sits the $4\times 4$ matrix
$B$ described above.

\noindent
{\bf The case $d=3$:}\\
The determinant of the form that appears in the following Theorem
is $2 \cdot 11 \cdot 6$, which is not a square,
hence there are no Steinhaus sets in dimension $3$.
\begin{theorem}
\label{th:3}
For each $x, y, z \in \ZZ$ the number
$$
Q(x,y,z) = 2x^2 + 11y^2 + 6z^2
$$
is a sum of three integer squares.
\end{theorem}

\noindent{\bf Proof.}
Suppose this is false and that there are $(x_0, y_0, z_0) \neq (0, 0, 0)$
and
\begin{itemize}
\item[(a)]
$Q(x_0, y_0, z_0)$ is not a sum of three squares, and
\item[(b)]
$x_0^2 + y_0^2 + z_0^2$ is minimal.
\end{itemize}
From (a), and the well known characterization of those natural numbers
that cannot be written as a sum of three squares, we have that
$$
Q(x_0, y_0, z_0) = 4^\nu(8k+7),\ \ \nu\ge0, k\ge0.
$$
If all $x_0, y_0, z_0$ are even, we have $\nu\ge 1$, and, setting
$x_0 = 2 x_1$, $y_0 = 2y_1$ and $z_0 = 2z_1$,
we obtain that $Q(x_1, y_1, z_1)$ is not a sum of three squares, which
contradicts the minimality of the initial triple $(x_0, y_0, z_0)$.
We conclude that at least one of $x_0, y_0, z_0$ is odd.

\noindent{\bf Case No 1:} $\nu=0$.

Then $Q(x_0, y_0, z_0) = 7 \bmod 8$.
But the quadratic residues mod $8$ are $0$, $1$ and $4$, and one
checks by examining all the possibilities that $Q$ is never $7 \bmod 8$.

\noindent{\bf Case No 2:} $\nu=1$.

Then $Q(x_0, y_0, z_0) = 32k + 28$.
Hence $y_0$ is even, say $y_0 = 2 y_1$.
We get
$$
x_0^2 + 22y_1^2 + 3 z_0^2 = 16k + 14,
$$
from which we conclude that $x_0$ and $z_0$ are odd,
$x_0 = 2x_1 + 1$, $z_0 = 2z_1 + 1$.
Substitution gives
\begin{eqnarray*}
4x_1^2 + 4x_1 + 1 + 22y_1^2 + 12 z_1^2 + 12 z_1 + 3 &=& 16k + 14\\
2x_1(x_1+1) + 11y_1^2 + 6z_1(z_1+1) + 2 &=& 8k+7\\
2x_1(x_1+1) + 11y_1^2 + 6z_1(z_1+1) &=& 5 \bmod 8.
\end{eqnarray*}
But $\xi^2+\xi = 0$ or $2$ or $4$ or $6 \bmod 8$, for all $\xi$,
hence, by applying this to the first and last term in the above sum,
and checking all possibilities we get a contradiction.

\noindent{\bf Case No 3:} $\nu\ge 2$.

As in Case No 2: $y_0 = 2y_1$, $z_0 = 2z_1+1$, $x_0 = 2x_1+1$.
Hence
$$
2x_1(x_1+1) + 11y_1^2 + 6z_1(z_1+1) + 2 = 4^{\nu-1}(8k+7),\ \ \nu-1\ge1.
$$ 
So $y_1$ is even, $y_1 = 2 y_2$, which gives
$$
x_1(x_1+1) + 22y_2^2 + 3 z_1(z_1+1) + 1 = 2\cdot 4^{\nu-2}(8k+7),
$$
a contradiction as the left hand side is odd while the right hand side
is even.
\Qed
We point out here that the actual quadratic form was only found by a
semi-automated computer search.
See \cite{mauldin-yingst} for a more systematic study of the method.

It is also shown in \cite{kolountzakis-papadimitrakis:quad}
that the method shown above cannot be applied
in dimension 2 to show the non-existence of
measurable sets with the Steinhaus property.
\begin{theorem}\label{th:2}
{\bf (Kolountzakis and Papadimitrakis, 2002)}\\
Any positive-definite binary quadratic form whose values are always sums of two
integer squares must have a determinant which is the square of an integer.
\end{theorem}

\subsection{Multi-lattice tiles}

\subsubsection{A ``finite'' Steinhaus problem}

The Steinhaus question essentially asks if there is a set
in the plane which is {\em simultaneously} a translational tile
for each translation set in the collection
$$
\Set{R_\theta \ZZ^2:\ 0\le\theta<2\pi},\ \ \ (\mbox{$R_\theta$ is rotation by $\theta$}).
$$
Restricting ourselves to the measurable case again it is easy to see,
using, for example, the Fourier method, that it is sufficient for a set to
be a tile for a countable dense (in the obvious sense) subset of these lattices (groups) in
order to be a tile for all of them.

The problem only becomes significantly different if one restricts oneself
to a {\em finite} collection of lattices
$$
\Lambda_0,\ldots,\Lambda_n \subseteq \RR^d,
$$
all of the same volume, say volume $1$, and asks for
a measurable subset of $\RR^d$ which tiles with all of them.
It turns out \cite{kolountzakis:multi} that this is generically feasible and we give here
a construction.
\begin{theorem}\label{th:multi-main}
{\bf (Kolountzakis, 1997)}\\
If the lattices $\Lambda_0,\ldots,\Lambda_n \subset \RR^d$
all have the same volume and if the sum of their dual lattices
$$
\Lambda_0^* + \cdots + \Lambda_n^*
$$
is a direct sum (i.e.\ there are no non-trivial relations
$\lambda_0+\cdots+\lambda_n = 0$ with $\lambda_i \in \Lambda_i^*$)
then they possess a Borel measurable common tile
(which is generally unbounded).
\end{theorem}

\noindent{\bf Proof.}
The common tile $\Omega \subseteq \RR^d$ that we construct
is a countable union of disjoint closed polyhedra (in fact, rectangles).

\begin{definition}
{\rm (Property A)}\\
We shall say that a collection of lattices
$\Lambda_0,\ldots,\Lambda_n \subset \RR^d$
has Property A if for each $\epsilon>0$ and for each
$x_0,\ldots,x_n \in \RR^d$ there exist
$\lambda_0 \in \Lambda_0, \ldots, \lambda_n \in \Lambda_n$,
with $\Abs{\lambda_j}$ arbitrarily large,
such that
\begin{equation}\label{eq:property-A}
\Abs{x_i - \lambda_i - (x_j - \lambda_j)} \le \epsilon,
\ \ \mbox{for all $i,j=0,\ldots,n$}.
\end{equation}
\end{definition}

That is, we can get any collection of points $x_0,\ldots,x_n \in \RR^d$
arbitrarily close to each other by translating $x_i$ by some
$\lambda_i\in\Lambda_i$, $i=0,\ldots,n$. 

We first show that if the given collection of lattices
has Property A then it has a common tile. At the end of the proof we
indicate why it is precisely the collections of lattices with their duals
having a direct sum that have Property A.

The letter $C$ will stand in this section for a positive constant
that may not depend on the parameter $K\to\infty$ and this
constant is not necessarily the same in all its occurences.

The lattices $\Lambda_j$, $j=0,\ldots,n$, are given by
\beql{lattice-dfn}
\Lambda_j = A_j \ZZ^d,\ \ \det A_j = 1.
\eeq
Let $D_j$ be the standard tile for the lattice $\Lambda_j$, i.e.,
\beql{std-FD}
D_j = A_j [0,1)^d,
\eeq
which is a parallelepiped of volume $1$.

Let $\Omega_0 = \emptyset$. In the end we shall have
$$
\Omega = \bigcup_{k=1}^\infty \Omega_k,
$$
where the $K$-th approximation
$$
A_K = \bigcup_{k=1}^K \Omega_k
$$
has measure $\mu(A_K) \to 1$, as $K \to \infty$,
and for each $j=0,\ldots,n$ almost all cosets
$x+\Lambda_j$ have no more than one point in $A_K$.
It follows that $\Omega$ contains exactly one element from
almost all the cosets of $\Lambda_j$, for each $j=0,\ldots,n$,
and is therefore a common tile for
the collection $\Lambda_0,\ldots,\Lambda_n$.
Assume that we have already defined $\Omega_0,\ldots,\Omega_K$.
The set $\Omega_{K+1}$ will be defined as follows.
The ``projection'' $\pi_j:\RR^d \to D_j$ is
defined by the relation
$$
x - \pi_j(x) \in \Lambda_j.
$$
The ``leftover'' after stage $K$
is then defined by
\begin{equation}\label{eq:leftovers}
L_j^{(K)} = D_j \setminus \pi_j(A_K),\ \ \mbox{for $j=0,\ldots,n$}.
\end{equation}
We have to ensure that $\mu\left(L_j^{(K)}\right) \to 0$, as $K \to \infty$.

\begin{figure}[h] \centering
\input{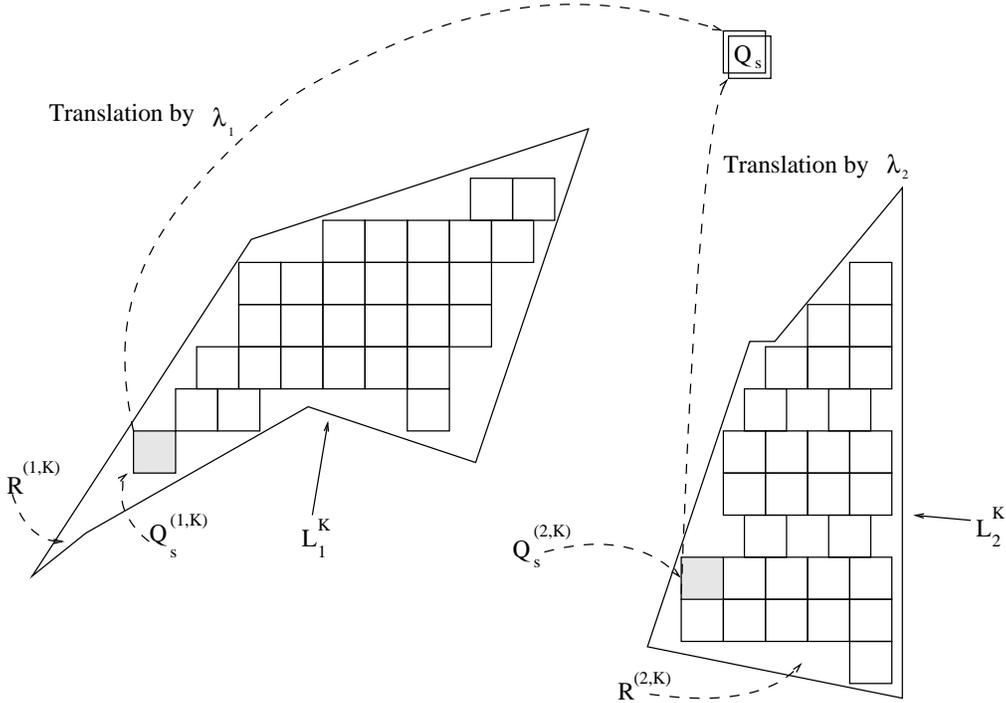}
\caption{Construction of the common tile for two lattices, $d=2$}
\ \\ \hrule
\label{fig:multi-construction}
\end{figure}

Our construction will guarantee that
each of the leftovers $L_j^{(K)}$ consists of a finite collection
of polyhedra.
Choose $\epsilon>0$ to be so small so as to be able
to write
\begin{equation}\label{eq:cut-up}
L_j^{(K)} = \left(\bigcup_{s=1}^{S} Q_s^{(j,K)}\right) ~\cup~ R^{(j,K)},
 \ \ \ (j=0,\ldots,n)
\end{equation}
where the $Q_s^{(j,K)}$, $s=1,\ldots,S=S(K)$, are axis-aligned, closed cubes
with disjoint interiors of side $\epsilon$, and
\begin{equation}\label{eq:R}
\mu\left(R^{(j,K)}\right) \le {1\over K}.
\end{equation}
Notice that the same number $S = S(K)$ of cubes is used independently of $j$.
(The conctruction is shown for two lattices in Figure \ref{fig:multi-construction}
in dimension $d=2$.)

For each $s=1,\ldots,S$, let $c_s^{(j,K)}$ be the center of
the cube $Q_s^{(j,K)}$ and, using Property A,
define $\lambda_s^{(j,K)} \in \Lambda_j$ to be such that all
$$
c_s^{(j,K)} - \lambda_s^{(j,K)},\ \ j=0,\ldots,n,
$$
are at most $\Ds {\epsilon \over K}$ apart.
The $\lambda_s^{(j,K)}$ are also taken large enough so that,
for fixed $j$,
no two translated cubes $Q_s^{(j,K)} - \lambda_s^{(j,K)}$ overlap.

Consider then the intersection of the $n+1$ translated cubes
\begin{equation}\label{eq:intersect}
\widetilde{Q}_s^{(K)} =
  \bigcap_{j=0}^n \left( Q_s^{(j,K)} - \lambda_s^{(j,K)} \right)
\end{equation}
and notice that
\begin{equation}\label{eq:interect-meas}
\mu(\widetilde{Q}_s^{(K)}) \ge \epsilon^d - C {\epsilon^d \over K}.
\end{equation}
Define
$$
\Omega_{K+1} = \bigcup_{s=1}^S \widetilde{Q}_s^{(K)}.
$$
We have $L_j^{(K+1)} = L_j^{(K)} \setminus \pi_j(\Omega_{K+1})$ and
$$
\mu\left( L_j^{(K)}\right) \to 0,
$$
as $K \to \infty$.
This is so because $L_j^{(K)} \setminus \pi_j(\Omega_{K+1})$
consists of the sets $R^{(j,K)}$, $j=0,\ldots,n$, which
have total measure $\Ds \le {n+1\over K}$
plus a set of measure $\Ds C {\epsilon^d \over K}$
for each $s=1,\ldots,S$, which amounts to no more
than $\Ds {C \over K}$ of measure, as clearly
$\epsilon^d S \le 1$.
\Qed

\begin{open-problem}
Can two lattices in generic position have a {\em bounded} measurable
common tile?
\end{open-problem}

\subsubsection{Multi-lattice tiles: an application to Weyl-Heisenberg bases}

\begin{definition}
{\bf (Gabor or Weyl-Heisenberg bases)}\\
A {\em Gabor} (or {\em Weyl-Heisenberg})
basis of $\RR^d$ is a function $g\in L^2(\RR^d)$, together
with two lattices $K=A\ZZ^d$ (the {\em translation lattice})
and $L=B\ZZ^d$ (the {\em modulation lattice}) such that the collection
\beql{gabor}
\Set{g(x-\kappa)e^{-2\pi i \lambda x},\ \ \kappa\in K,\ \lambda\in L},
\eeq
is an orthonormal basis of $L^2(\RR^d)$.
\end{definition}

It had been known for some time (see the introduction and references
in \cite{han-wang}) that if there is a Weyl-Heisenberg basis for
the lattices $K$ and $L$ then it must be true that
\begin{equation}\label{gabor-density-condition}
\dens K \cdot \dens L = 1.
\end{equation}
Apart from dimension $1$ though, the converse had not been known
until Han and Wang \cite{han-wang} used the idea of multi-lattice
tiles to prove that whenever \eqref{gabor-density-condition} holds
then there is a $g$ such that collection \eqref{gabor} is an orthonormal
basis of $L^2(\RR^d)$.

Han and Wang \cite{han-wang} first proved that the genericity condition
described in Theorem \ref{th:multi-main} is not necessary when the number
of lattices is two.
\begin{theorem}\label{th:han-wang-two}
{\bf (Han and Wang, 2001)}\\
Whenever the lattices $\Lambda_0$ and $\Lambda_1$ in $\RR^d$ have the same
volume then there exists a measuarable set $E \subseteq \RR^d$
which tiles with both of them
\end{theorem}
Thus, for two lattices of the same volume there is
{\em always} a measurable common tile. This is {\em not} true for three
or more lattices without some condition,
as the following result \cite{kolountzakis:multi} shows.

\begin{theorem}\label{th:commensurable-lattices}
{\bf (Kolountzakis, 1997)}\\
There are three lattices in $\RR^2$ which have the same volume
and do not admit a common tile.
\end{theorem}

\noindent{\bf Proof.}
Let
$$
\Lambda_0 = (2\ZZ)\times\ZZ,\ \Lambda_1 = \ZZ\times(2\ZZ),
\ \mbox{and}\ \Lambda_2 = \Set{(k,l)\in\ZZ^2:\ k=l \bmod 2}.
$$
It is easy to see that
$$
\ZZ^2 = \sum_{i=0}^2 \Lambda_i = \bigcup_{i=0}^2 \Lambda_i.
$$
Suppose now that $\Omega\subset\RR^2$ is such that for all
$x\in\RR^2$, outside a set $E$ of measure $0$, we have that
$x+\Lambda_i$ contains exactly one point of $\Omega$, for all
$i=0,1,2$.
(We do not assume that $\Omega$ is measurable.)
It follows that for almost all $x\in\RR^2$ (with an exceptional set
perhaps different from $E$) we have
$$
\Abs{(x+\ZZ^2)\cap\Omega}=2\mbox{\ and\ }
\Abs{(x+\Lambda_i)\cap\Omega}=1,\ i=0,1,2.
$$
Indeed, $\ZZ^2$ is the disjoint union of $\Lambda_0$ and
$\Lambda_0+(1,0)$ and so are all its translates.
We define the set
$$
E' = E \cup (E-(1,0)),
$$
which is clearly still a null set.
Then, for $x \notin E'$ the set $x+\ZZ^2$ contains exactly
two points of $\Omega$, since the two disjoint copies of $\Lambda_0$
therein both contain exactly one $\Omega$-point.

By translating $\Omega$ we may assume that this holds for $x=0$.
Let then
$$
\Set{z,w} = \ZZ^2\cap\Omega.
$$
It follows that $z-w \in \ZZ^2$ and, since $\ZZ^2 = \bigcup_{j=1}^3 \Lambda_j$,
$z-w$ belongs to some $\Lambda_j$.
But then the $\Omega$-points $z$ and $w$ belong to the same $\Lambda_j$-coset,
a contradiction.
Hence the $\Lambda_i$ have no common tile in $\RR^2$ in a strong
sense.
\Qed

We continue now with proof of Han and Wang \cite{han-wang}
that \eqref{gabor-density-condition} suffices for the existence of a fucntion
$g$ such that the collection \eqref{gabor} is a Weyl-Heisenberg basis.
Suppose then that \eqref{gabor-density-condition} holds. It follows
that the lattices $K$ and $L^*$ have the same volume.
Hence, by Theorem \ref{th:han-wang-two},
there is a common tile $E\subseteq\RR^d$ for $K$ and $L^*$.
Let
$$
g = \chi_E.
$$
For any $f\in L^2(\RR^d)$ write then
$$
f(x) = \sum_{\kappa \in K} f_\kappa(x) := \sum_{\kappa\in K} g(x-\kappa) f(x)
$$
which is an orthogonal decomposition precisely because $E$ is a $K$-tile.
For each $\kappa$, $f_\kappa(x)$ is a function on $E+\kappa$ which is a
$L^*$-tile.
But if a set $\Omega$ tiles with a lattice $L^*$ then the collection
$$
\Set{\exp^{2\pi i \inner{\lambda}{x}}:\ \lambda \in L}
$$
is an orthogonal basis for $L^2(\Omega)$
(this is merely multi-dimensional Fourier Series plus a change of variable,
but see also Theorem \ref{th:fuglede-lattice} below).
For $\Omega = E+\kappa$ we therefore obtain that
$$
f_\kappa(x) = \sum_{\lambda\in L} \inner{f_\kappa}{e^{2\pi i\inner{\lambda}{x}}}
 e^{2\pi i\lambda x}\ \ \ (x \in E+\kappa)
$$
is an orthogonal decomposition and so is then
$$
f(x) = \sum_{\kappa\in K, \lambda\in L} \inner{f}{g(x-\kappa)e^{2\pi i\inner{\lambda}{x}}}
 g(x-\kappa)e^{2\pi i\lambda x},
$$
as as we had to show.

\subsection{The support of ``soft'' multi-lattice tiles}

Fix the dimension $d$ and
take {\em any} finite collection of lattices $\Lambda_1,\ldots,\Lambda_N$.
Then the function
\beql{convolution-of-tiles}
f = \chi_{D_1}*\cdots*\chi_{D_N},
\eeq
where $D_j$ is any tile for $\Lambda_j$, tiles with the given lattices,
as one can see directly from the definition of tiling (if $f+\Lambda$ is a
tiling then so is $f*g + \Lambda$, even for non-lattice $\Lambda$).

For this particular $f$ (and whatever choice of $D_j$) we have
$$
{\rm diam}\,\supp f \ge C N,
$$
with a constant that dependes only on $d$. This is easy to see as at least
$1/d$ of the sets $D_j$ will be ``long'' along the same one of the $d$ coordinate
axes and the convolution of all of them will therefore also be long along that
axis.

If one chooses appropriate parallelograms for the $D_j$'s one gets more
or less the best known (to me at least) construction as regards the diameter
of the common tile of the collection $\Lambda_1,\ldots,\Lambda_N$, where, now,
we do not insist that the tile be an indicator function,
but rather any integrable function. One can in this manner get a tile whose
support has diameter $\sim N$.

It is not obvious at all that this size has to grow as a function of $N$.
In fact, the following theorem \cite{kolountzakis-wolff}, which provides
a lower bound for the diameter of the support of a common tile,
is the only one of its kind, uses 
(multivariable) entire function theory (some times ineffective in such matters) and is still
far from the best known upper bound ($\sim N$).
\begin{theorem}
\label{th:kw-lower-bounds}
{\bf (Kolountzakis and Wolff, 1997)}\\
Suppose that $\Lambda_1, \ldots, \Lambda_N$ are unimodular lattices
in $\RR^d$ with $\Lambda_i \cap \Lambda_j = \Set{0}$ for all $i\neq j$.
Suppose also that the non-zero $f\in L^1(\RR^d)$ is a common tile for the $\Lambda_j$.
Then
$$
{\rm diam}\,\supp f \ge C_d N^{1/d}.
$$
\end{theorem}

\noindent{\bf Proof.}
All constants below may depend only on the dimension $d$.
We note that $\Lambda_1 \cap \Lambda_2 = \Set{0}$ implies that
the lattice $\Lambda_1^*$ is uniformly distributed mod $\Lambda_2^*$.
This can be proved using Weyl's lemma--see for example
\cite{kolountzakis:multi}.

We shall make use of a theorem of  Ronkin \cite{ronkin}
and  Berndtsson \cite{berndtsson} which concerns the zero set
on the real plane of an entire function of several complex variables
which is of exponential type. We formulate it as a lemma:
\begin{lemma}\label{ronkin-berndtsson}
{\bf (Ronkin 1972, Berndtsson 1978)}\\
Assume that $E\subset\RR^d$ is a countable set with any two
points having distance at least $h$ and let
$$
d_E = \limsup_{r\to\infty} {\Abs{E \cap D(0,r)}\over \Abs{D(0,r)}}
$$
be its upper density (see Definition \ref{def:density}).
Assume that $g:\CC^d\to\CC$ is an entire function vanishing
 on $E$
which is of exponential type
$$
\sigma < A(d) h^{d-1} d_E.
$$
Then $g$ is identically $0$.
(Here $A(d)$ is an explicit function of the dimension $d$.)
\end{lemma}
When $d=1$ this is  classical and follows  from Jensen's formula.

Assume that $f:\RR^d\to\CC$ is as in Theorem 3 and write
$$
\alpha = \mbox{diam supp}\, f
$$
We may assume  that supp$f$ is contained in a disc of radius
$\lesssim\alpha$ centered at the origin, since the assumptions
 are unaffected by a translation of coordinates. Then
  $\ft{f}$ can be extended to $\CC^d$
as an entire function of exponential type $C\alpha$, in fact
$$
\Abs{\ft{f}(x+iy)} \le C_f e^{C \alpha\Abs{y}},
\ \ \mbox{for $x+iy \in \CC^d$}.
$$
Furthermore, since $f$ tiles with all $\Lambda_j$, it follows that
$\ft{f}$ vanishes on
$$
{\cal Z} = \bigcup_{i=1}^n\Lambda_i^*\setminus\Set{0}.
$$
Observe that, since every lattice $\Lambda_i^*$ is
uniformly distributed mod every $\Lambda_j^*$, $j \neq i$,
the density of points in each $\Lambda_i^*$ which are also
in some $\Lambda_j^*$ is $0$ and therefore the density of the set ${\cal Z}$
is equal to $n$.

In order to use Lemma \ref{ronkin-berndtsson} we have to select
a large (in terms of upper density), well-separated subset
of ${\cal Z}$.
Notice first that we can assume that for each $i$ all points of
$\Lambda_i^*$ are at least distance $ n^{-\frac{1}{d}}$ apart.
For if $u,v \in \Lambda_i^*$ have $\Abs{u-v} <  n^{-\frac{1}{d}}$ then,
for a suitable constant $c$, the one-dimensional version of Lemma \ref{ronkin-berndtsson}
implies that the function
$\ft{f}$ on the subspace $E = \CC(u-v)$ cannot be of exponential
type $\leq   cn^{\frac{1}{d}}$.
Indeed, $\ft{f}$ would have too many zeros on that subspace, namely all
multiples of $u-v$, which all belong to $\Lambda_i^*$.
Note also that $\ft{f}$ does not vanish identically on this
subspace. But $\ft{f}$ restricted to $E$
is the Fourier Transform of $f_E:E\to\CC$ defined by
$f_E(x) = \int_{x+E^\perp} f(y)\,dy$ (here $E^\perp$ is the orthogonal complement
of $E\cap\RR^n$ in $\RR^n$).
Hence $\alpha \ge \mbox{diam supp}\, f_E \ge C n^{\frac{1}{d}}$,
which is what we want to conclude about $\alpha$.

Suppose now that we want to extract a subset of ${\cal Z}$ whose elements
are at least $h$ distance apart, for some $h>0$ to be determined
later.
We shall say that point $x$ of lattice $\Lambda_i^*$ is {\em killed}
by point $y$ of lattice $\Lambda_j^*$ if $\Abs{x-y}<h$.
Then, we define the subset ${\cal Z}'$ of ${\cal Z}$ as those points of ${\cal Z}$
which are not killed by any point of the {\em other} lattices.
This set clearly has all its points at distance at least $h$ apart,
provided that
\beql{h-small}
h \le {1\over2} \min_{u,v \in \Lambda_i^*} \Abs{u-v} \le C
 n^{-\frac{1}{d}},
\eeq
so that no point of a lattice may kill a point of the same lattice.
Let us see how many points of $\Lambda_2^*$ are killed by
some point of $\Lambda_1^*$.
We use the uniform distribution of $\Lambda_2^*$ mod $\Lambda_1^*$.

Fix a fundamental parallelepiped $D_1$ of $\Lambda_1^*$.
It is clear that only a fraction
$\rho(h) \le C h^d$ of $D_1=\RR^d/\Lambda_1^*$
has distance from $0$ that is less than $h$ (this distance is measured
on the {\em torus} $D_1$).
As $\Lambda_2^*$ is uniformly distributed mod $\Lambda_1^*$
the subset of points of $\Lambda_2^*$ which are killed
by some point of $\Lambda_1^*$ has density $\rho(h)$.
Hence the density of those points of $\Lambda_2^*$ that are
killed by {\em any} other lattice is at most
$(n-1)\rho(h) \le C h^d n$.
We deduce that the density of ${\cal Z}'$ is
at least $(1-C n h^d) n$.
We now choose $h=c n^{-\frac{1}{d}}$, for a sufficiently small
 constant $c$,
to ensure that the density of ${\cal Z}'$ is at least
$C n$.
Applying  Lemma \ref{ronkin-berndtsson} with $g=\ft{f}$
and $E = {\cal Z}'$ we get
$$
\alpha \ge C A h^{d-1} n \ge C n^{\frac{1}{d}}.
$$
\Qed

\begin{open-problem}
Bridge the gap between Theorem \ref{th:kw-lower-bounds}
and the upper bound $\sim N$.
\end{open-problem}

\pagebreak
\section{Lecture 3: The Fuglede Conjecture}

\subsection{Spectral sets and tiling}

Let us write $e_\lambda(x) = \exp 2\pi i\inner{\lambda}{x}$.
\begin{definition}\label{def:spectral-set}
{\bf (Spectral sets)}\\
Suppose that $\Omega$ is a bounded open set of measure $1$.
We call $\Omega$ {\em spectral} if $L^2(\Omega)$
has an orthonormal basis
$$
E_\Lambda=\Set{e_\lambda:\ \lambda\in\Lambda}
$$
of exponentials.
The set $\Lambda$ is then called a {\em spectrum} for $\Omega$.
\end{definition}
(We only restrict ourselves to sets of measure $1$ to make our life simpler.)

The inner product and norm on $L^2(\Omega)$ are
$$
\inner{f}{g}_\Omega = \int_\Omega f \overline{g},
\ \mbox{ and }\
\Norm{f}_\Omega^2 = \int_\Omega \Abs{f}^2.
$$
We have
$$
\inner{e_\lambda}{e_x}_\Omega = \ft{\chi_\Omega}(x-\lambda).
$$
which gives
$$
\mbox{$E(\Lambda)$ is orthogonal} \Leftrightarrow
 \forall \lambda,\mu\in\Lambda, \lambda\neq\mu:
   \ \ft{\chi_\Omega}(\lambda-\mu) = 0
$$
For $E(\Lambda)$ to be complete as well we must in addition have
\beql{parseval}
\forall f\in L^2(\Omega):\ \ 
    \Norm{f}_2^2 =  \sum_{\lambda\in\Lambda} \Abs{\inner{f}{e_\lambda}}^2.
\eeq
It is sufficient to have \eqref{parseval} for $f(t) = e_x(t)$, $x\in\RR^d$,
since then we have it in the closed linear span of these functions, which is all
of $L^2(\Omega)$.

\begin{figure}[h] \centering
\input{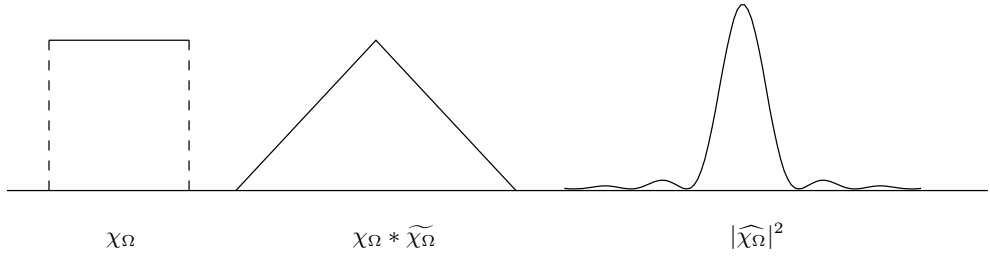}
\caption{The functions $\chi_\Omega$, $\chi_\Omega*\widetilde{\chi_\Omega}$
and $\Abs{\ft{\chi_\Omega}}^2$, when $\Omega$ is an interval.
The last two functions are a Fourier Transform pair.}
\ \\ \hrule
\label{fig:fuglede-interval}
\end{figure}

An equivalent reformulation for $\Lambda$ to be a spectrum of $\Omega$
is therefore the following, which we state as a theorem.
\begin{theorem}\label{th:spectral-reformulation}
The set $\Lambda$ is a spectrum of $\Omega$
if and only if $\sum_{\lambda\in\Lambda} \Abs{\ft{\chi_\Omega}(x-\lambda)}^2 = 1$,
 for almost every $x\in\RR^d$.
In tiling language
$$
\mbox{$\Lambda$ is a spectrum of $\Omega$} \ \Leftrightarrow\ 
\mbox{$\Abs{\ft{\chi_\Omega}}^2 + \Lambda = \RR^d$}
$$
\end{theorem}
The relevant functions are shown in Figure \ref{fig:fuglede-interval}, for
the case of $\Omega$ being an interval.

It follows from Theorem \ref{th:spectral-reformulation} that the spectrum
$\Lambda$ of domain $\Omega$, if it exists, has all the nice properties of
tiling sets.
In particular, $\Lambda$ has uniform density equal to $1$ and its points
are $\epsilon$-separated for some $\epsilon>0$. 

We can now state Fuglede's Conjecture \cite{fuglede}
\begin{conjecture}{\bf (Fuglede 1974)}\\
Let $\Omega\subseteq\RR^d$ be a bounded, open domain of measure $1$.
Then $\Omega$ is spectral
if and only if it can tile space by translation.
\end{conjecture}
We should emphasize here that no relation is claimed in the conjecture
between the spectrum of $\Omega$ and the set of translations with which $\Omega$ tiles.

\begin{remark}
By the preceding discussion Fuglede's Conjecture states that $\Omega$
is a tile if and only if $\Abs{\ft{\chi_\Omega}}^2$ is a tile
(both tilings are at level $1$).
\end{remark}

Despite a lot of work that has been done in the last 5-6 years
the conjecture remains open in all dimensions and in both directions.
One easy and important case though is given by the following \cite{fuglede}.
\begin{theorem}\label{th:fuglede-lattice}
{\bf (Fuglede, 1974)}\\
Suppose $\Omega\subseteq\RR^d$ is a bounded open domain of measure 1 and
$\Lambda\subseteq\RR^d$ a lattice of density 1.
Then $\Omega + \Lambda = \RR^d$ if and only if
$\Lambda^*$ (the dual lattice) is a spectrum of $\Omega$.
\end{theorem}

\noindent{\bf Proof.}
As remarked above, $\Lambda^*$ is a spectrum of $\Omega$
if and only if (see \S\ref{sec:lattice})
$$
\Abs{\ft{\chi_\Omega}}^2 + \Lambda^* = \RR^d,
$$
which is in turn equivalent to the Fourier Transform of the function $\Abs{\ft{\chi_\Omega}}^2$
vanishing on the dual lattice of $\Lambda^*$ except at $0$.
That is the function $f = \chi_\Omega * \widetilde{\chi_\Omega}$
vanishes on $\Lambda\setminus\Set{0}$.
But $f$ is non-zero exactly on $\Omega-\Omega$, hence the above vanishing
is equivalent to
$$
(\Omega-\Omega) \cap \Lambda = \Set{0}
$$
which means precisely that the copies $\Omega+\lambda$, $\lambda\in\Lambda$,
are non-overlapping. But $\Abs{\Omega}=1$ and $\dens\Lambda=1$, hence
the above packing is indeed a tiling.
The argument is completely reversible.
\Qed

\subsection{Implications of the Brunn-Minkowski inequality for
convex tiles and spectral bodies}

Let us recall a simple case of the {\bf Brunn-Minkowski Inequality}
(see e.g.\ \cite{schneider}).
\begin{quotation}
For a convex body $K$ we always have
$$
\Abs{K-K} \ge 2^d\Abs{K}.
$$
We have equality above exactly when $K$ is symmetric,
in which case $K-K = 2K$.
\end{quotation}
Using the Brunn-Minkowski inequality one can show:
\begin{theorem}\label{th:tiles-are-symmetric}
{\rm ({\bf Minkowski, ca.\ 1900})}\\
If $\Omega$ is a convex translational tile then it is symmetric.
\end{theorem}

\noindent{\bf Proof.}
Suppose $K$ is convex and $K+\Lambda = \RR^d$.
By the packing condition (non-overlaping of translates) only we get
$$
(K-K) \cap (\Lambda-\Lambda) = \Set{0}.
$$
Define the convex set $L = {1\over 2}(K-K)$.
One easily sees that $L-L = K-K$, so that
$$
L + \Lambda \le \RR^d,
$$
is a packing. But this implies (see Lemma \ref{lm:packing-implies-density-bound})
$$
\Abs{L} \le 1,
$$
and by the equality case in the Brunn-Minkowski inequality $K$ is symmetric.
\Qed

The following theorem \cite{kolountzakis:nonsym}
is also a consequence of the Brunn-Minkowski
inequality. 
\begin{theorem}\label{th:nonsym}
{\bf (Kolountzakis, 2000)}\\
If $\Omega$ is convex and spectral then it is symmetric.
\end{theorem}
This result is of course in agreement with the Fuglede Conjecture as this
would be false if there were any non-symmetric convex spectral domains.
We prove Theorem \ref{th:nonsym} in \S \ref{sec:fourier} and \S \ref{sec:nonsym-proof}
below.

\subsubsection{Fourier-analytic conditions for tiling}
\label{sec:fourier}

When studying tiling by the function $\Abs{\ft{\chi_\Omega}}^2$ Theorem
\ref{th:tiling-implies-condition} is not applicable since the
Fourier Transform of the function, namely $\chi_\Omega*\widetilde{\chi_\Omega}$,
is never smooth.
However, the positivity of the function and its Fourier Transform as well as the
compact support of the Fourier Transform compensate for this lack of smoothness and allow
us to prove the following result \cite{kolountzakis:nonsym}.
\begin{theorem}\label{th:tiling-implies-supp}
{\bf (Kolountzakis, 2000)}\\
Suppose that $f\ge 0$ is not identically $0$, that $f \in L^1(\RR^d)$,
$\widehat{f}\ge 0$ has compact support and $\Lambda\subset\RR^d$.
If $f+\Lambda$ is a tiling then
\beql{sp-cond-1}
\supp \ft{\delta_\Lambda} \subseteq \Set{x\in\RR^d:\ \ft{f}(x) = 0} \cup \Set{0}.
\eeq
\end{theorem}

\noindent{\bf Proof.}
Assume that $f+\Lambda = w\RR^d$ and let
$$
K = \Set{\widehat{f} = 0} \cup \Set{0}.
$$
We have to show that
$$
\widehat{\delta_\Lambda}(\phi) = 0,
\ \ \forall \phi\in C_c^\infty(K^c).
$$
Since $\ft{\delta_\Lambda}(\phi) = \delta_\Lambda(\ft\phi)$
this is equivalent to $\sum_{\lambda\in\Lambda}\widehat\phi(\lambda) = 0$, for
each such $\phi$.
Notice that $h = \phi / \widehat{f}$ is a continuous function, but
not necessarily smooth. We shall need that $\widehat h \in L^1$.
This is a consequence of a well-known theorem of Wiener \cite[Ch.\ 11]{rudin}.
We denote by $\TT^d = \RR^d/\ZZ^d$ the $d$-dimensional torus.
\begin{theorem}\label{th:wiener}
{\bf (Wiener)}\\
If $g \in C(\TT^d)$ has an absolutely convergent Fourier series
$$
g(x) = \sum_{n\in\ZZ^d} \widehat g(n) e^{2\pi i \Inner{n}{x}},
\ \ \ \widehat g \in \ell^1(\ZZ^d),
$$
and if $g$ does not vanish anywhere on $\TT^d$ then
$1/g$ also has an absolutely convergent Fourier series.
\end{theorem}
Assume that
$$
\supp\phi,\ \supp\widehat{f} \subseteq \left(-{L\over2},{L\over2}\right)^d.
$$
Define the function $F$ to be:\\
(i) periodic in $\RR^d$ with period lattice $(L\ZZ)^d$,\\
(ii) to agree with $\widehat{f}$ on $\supp\phi$,\\
(iii) to be non-zero everywhere and,\\
(iv) to have $\widehat{F} \in \ell^1(\ZZ^d)$, i.e.,
$$
\widehat F = \sum_{n\in\ZZ^d} \widehat F(n) \delta_{L^{-1} n},
$$
is a finite measure in $\RR^d$.

One way to define such an $F$ is as follows.
First, define the $(L\ZZ)^d$-periodic function $g\ge 0$
to be $\ft{f}$ periodically extended.
The Fourier coefficients of $g$ are $\ft{g}(n) = L^{-d} f(-n/L) \ge 0$.
Since $g, \ft g \ge 0$ and $g$ is continuous at $0$ it is
easy to prove that $\sum_{n\in\ZZ^d} \ft g (n) = g(0)$, and therefore
that $g$ has an absolutely convergent Fourier series.

Let $\epsilon$ be small enough to guarantee that $\ft{f}$
(and hence $g$) does not vanish on $(\supp\phi) + B_\epsilon(0)$.
Let $k$ be a smooth $(L\ZZ)^d$-periodic function which is equal to $1$
on $(\supp\phi)+(L\ZZ^d)$ and equal to $0$ off
$(\supp\phi + B_\epsilon(0))+(L\ZZ^d)$, and
satisfies $0\le k \le 1$ everywhere.
Finally, define
$$
F = k g + (1-k).
$$
Since both $k$ and $g$ have absolutely summable Fourier series and this
property is preserved under both sums and products, it follows that $F$ also
has an absolutely summable Fourier series. And by the nonnegativity of $g$ we
get that $F$ is never $0$, since $k=0$ on $\Set{\ft f=0}+(L\ZZ^d)$.

By Wiener's Theorem \ref{th:wiener}, $\widehat{F^{-1}} \in \ell^1(\ZZ^d)$, i.e., $\ft{F^{-
1}}$ is a finite
measure on $\RR^d$.
We now have that
$$
\FT{{\phi \over \widehat f}} = \widehat{\phi F^{-1}} =
  \widehat\phi * \widehat{F^{-1}} \in L^1(\RR^d).
$$
This justifies the interchange of the summation and integration below:
\begin{eqnarray*}
\sum_{\lambda\in\Lambda} \widehat\phi(\lambda) 
 &=& \sum_{\lambda\in\Lambda} \FT{{\phi\over\widehat{f}} \widehat{f}}
(\lambda) \\
 &=& \sum_{\lambda\in\Lambda} \FT{{\phi\over\widehat f}} *
\widehat{\widehat{f}}~(\lambda) \\
 &=& \sum_{\lambda\in\Lambda} \int_{\RR^d} \FT{{\phi\over\widehat f}}(y) f(y-
\lambda) ~dy \\
 &=& \int_{\RR^d} \FT{{\phi\over\widehat f}}(y)
  \sum_{\lambda\in\Lambda} f(y-\lambda) ~dy\\
 &=& w \int_{\RR^d} \FT{{\phi\over\widehat f}}(y)~dy\\
 &=& w{\phi\over\ft{f}}(0)\\
 &=& 0,
\end{eqnarray*}
as we had to show.
\Qed

For a set $A \subseteq \RR^d$ and $\delta > 0$ we write
$$
A_\delta = \Set{x\in\RR^d:\ \dist(x, A)<\delta}.
$$

\begin{figure}[h] \centering
\input{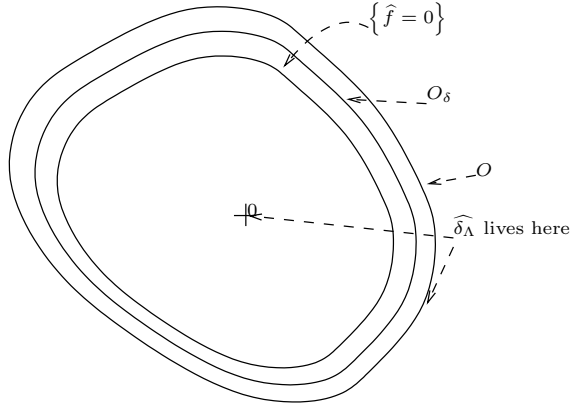}
\caption{The sets appearing in Theorem \ref{th:disjoint-supports}.
The sets $O, O_\delta, \Set{\ft{f}=0}$ all live {\em outside} the contours.}
\ \\ \hrule
\label{fig:zero-set-fuglede}
\end{figure}

We shall need the following partial converse to Theorem
\ref{th:tiling-implies-supp} (see Figure \ref{fig:zero-set-fuglede}
for the assumptions of Theorem \ref{th:disjoint-supports}).
\begin{theorem}\label{th:disjoint-supports}
Suppose that $f\in L^1(\RR^d)$,
and that $\Lambda \subset \RR^d$ has uniformly bounded density.
Suppose also that $O \subset \RR^d$ is open and
\beql{dijoint-supports}
\supp \ft{\delta_\Lambda} \setminus \Set{0} \ \subseteq\ O
\ \mbox{ and }\ O_\delta \subseteq\  \Set{\ft{f} = 0},
\eeq
for some $\delta>0$.
Then $f+\Lambda$ is a tiling at level $\ft{f}(0)\cdot\ft{\delta_\Lambda}(\Set{0})$.
\end{theorem}
The assumptions of Theorem \ref{th:disjoint-supports} ensure that the supports
of $\ft{\delta_\Lambda}$ (except at $0$) and $\ft f$ are well separated.
In other words $\ft f$ vanishes to infinite order on the support of $\ft{\delta_\Lambda}$.
This makes the formal implication
$$
\ft f \cdot \ft{\delta_\Lambda} = \ell \delta_0 \Longrightarrow
f*\delta_\Lambda = \ell
$$
correct.
\begin{remark}
By the assumptions of the theorem we know that $\ft{\delta_\Lambda}$
is supported only at $0$, in a neighborhood of the origin.
It follows from Theorem \ref{th:point-implies-measure} that
$\ft{\delta_\Lambda}$ is a measure in some neighborhood of the origin
so it makes sense to speak of $\ft{\delta_\Lambda}(\Set{0})$.
\end{remark}

\noindent{\bf Proof.}
Let $\psi:\RR^d\to\RR$ be smooth, have support in $B_1(0)$ and $\ft\psi(0)=1$ and for
$\epsilon>0$ define the approximate identity
$\psi_\epsilon(x) = \epsilon^{-d}\psi(x/\epsilon)$.
Let
$$
f_\epsilon = \ft{\psi_\epsilon} f,
$$
which has rapid decay.

First we show that $(\int f_\epsilon)^{-1}f_\epsilon + \Lambda$
is a tiling.
That is, we show that the convolution $f_\epsilon * \delta_\Lambda$
is a constant.
Let $\phi$ be any Schwartz function. Then
$$
f_\epsilon * \delta_\Lambda (\phi) =
\ft{f_\epsilon}\ft{\delta_\Lambda} (\ft\phi(-x)) =
\ft{\delta_\Lambda} (\ft\phi(-x) \ft{f_\epsilon}).
$$
The function $\ft\phi(-x) \ft{f_\epsilon}$ is a Schwartz
function whose support intersects $\supp\ft{\delta_\Lambda}$ only at $0$,
since, for small enough $\epsilon>0$,
$$
\supp \ft\phi \ft{f_\epsilon} \subseteq
\supp \ft{f_\epsilon} \subseteq
(\supp\ft{f})_\epsilon \subseteq O^c.
$$
Hence, for each Schwartz function $\phi$
$$
f_\epsilon * \delta_\Lambda (\phi) = 
	\ft\phi(0) \ft{f_\epsilon}(0) \ft{\delta_\Lambda}(\Set{0}),
$$
which implies
$$
f_\epsilon * \delta_\Lambda (x) = \ft{f_\epsilon}(0) \ft{\delta_\Lambda}(\Set{0}),
\ \ \mbox{a.e.($x$)}.
$$
We also have that $\sum_{\lambda\in\Lambda} \Abs{f(x-\lambda)}$ is finite
a.e.\ (see the remark following the definition of tiling), hence,
for almost every $x\in\RR^d$
$$
\sum_{\lambda\in\Lambda}\Abs{f(x-\lambda) - f_\epsilon(x-\lambda)} =
\sum_{\lambda\in\Lambda}\Abs{f(x-\lambda)}\cdot\Abs{1-\ft{\psi_\epsilon}(x-\lambda)},
$$
which tends to $0$ as $\epsilon\to 0$.
This proves
$$
\sum_{\lambda\in\Lambda} f(x-\lambda) =
	\ft{f}(0)\cdot\ft{\delta_\Lambda}(\Set{0}),\ \ \mbox{a.e.($x$)}.
$$
\Qed

\subsubsection{Convex spectral bodies must be symmetric}
\label{sec:nonsym-proof}

\noindent{\bf Proof of Theorem \ref{th:nonsym}:}
Write $K = \Omega - \Omega$, which is a symmetric, open convex set.
Assume that $(\Omega, \Lambda)$ is a spectral pair.
We can clearly assume that $0 \in \Lambda$.
It follows that $\Abs{\ft{\chi_\Omega}}^2 + \Lambda$ is a tiling and hence that 
$\Lambda$ has uniformly bounded density, has density equal to $1$ and
$\ft{\delta_\Lambda}(\Set{0}) = 1.$

By Theorem \ref{th:tiling-implies-supp} (with
$f = \Abs{\ft{\chi_\Omega}}^2,\ \ \ft{f} =
	\chi_\Omega * \widetilde{\chi_\Omega}(-x)$)
it follows that
$$
\supp \ft{\delta_\Lambda} \subseteq \Set{0} \cup K^c.
$$
Let $H = K/2$ and write
$$
f(x) = \chi_H*\widetilde{\chi_H}(x) = \int_{\RR^d} \chi_H(y) \chi_H(y-x)~dy.
$$
The function $f$ is supported in $\overline{K}$ and has nonnegative Fourier Transform
$$
\ft{f} = \Abs{\ft{\chi_H}}^2.
$$
We have
$$
\int_{\RR^d}\ft{f} = f(0) = \vol H
$$
and
$$
\ft{f}(0) = \int_{\RR^d} f = (\vol H)^2.
$$
By the Brunn-Minkowski inequality
for any convex body $\Omega$,
$$
\vol {1\over2}(\Omega-\Omega) \ge \vol \Omega,
$$
with equality only in the case of symmetric $\Omega$.
Since $\Omega$ has been assumed to be non-symmetric it follows that
$$
\vol H > 1.
$$
For
$$
1>\rho>\left({1\over \vol H}\right)^{1/d}
$$
consider
$$
g(x) = f(x/\rho)
$$
which is supported properly inside $K$, and has
$$
g(0) = f(0) = \vol H,\ \ \int_{\RR^d} g = \rho^d \int_{\RR^d} f = \rho^d (\vol H)^2.
$$
Since $\supp g$ is properly contained in $K$
Theorem \ref{th:disjoint-supports} implies that $\ft{g} + \Lambda$ is a tiling
at level $\int\ft{g} \cdot \dens\Lambda = \int\ft{g} = g(0) = \vol H$.
However, the value of $\ft{g}$ at $0$ is $\int g = \rho^d (\vol H)^2 > \vol H$,
and, since $\ft{g} \ge 0$ and $\ft{g}$ is continuous, this is a contradiction.
\Qed

\subsection{The spectra of the cube}

In this section we prove the following
\cite{iosevich-pedersen,lagarias-reeds-wang,kolountzakis:spectral}.
\begin{theorem}\label{th:spectra-of-the-cube}
{\bf (Iosevich and Pedersen, 1998, Lagarias, Reeds and Wang 1998, Kolountzakis 1999)}\\
Let $Q=(-1/2, 1/2)^d$ be the unit cube in $\RR^d$
and $\Lambda\subseteq\RR^d$. Then
$$
\mbox{$\Lambda$ is a spectrum of $Q$} \Leftrightarrow
\mbox{$Q+\Lambda=\RR^d$}.
$$
\end{theorem}
This had been proved earlier by Jorgensen and Pedersen \cite{jorgensen-pedersen}
for $d=3$.

\subsubsection{A lemma for two different tiles}

The following simple result is rather unexpected.
It is intuitively clear when $\Lambda$ is a periodic set but it is,
perhaps, suprising that it holds without any assumptions on
the set $\Lambda$.

\begin{lemma}\label{lm:tiling-from-packing}
If $f, g \ge 0$, $\int f(x) dx = \int g(x) dx = 1$ and 
both $f+\Lambda$ and $g+\Lambda$ are packings of $\RR^d$, then
$f+\Lambda$ is a tiling if and only if $g+\Lambda$ is a tiling.
\end{lemma}

\noindent{\bf Proof.}
We first show that, under the assumptions of the Theorem,
\beql{supp-tiling}
\mbox{$f+\Lambda$ tiles $-\supp g$ } \Longrightarrow
\mbox{ $g+\Lambda$ tiles $-\supp f$}.
\eeq
Indeed, if $f+\Lambda$ tiles $-\supp g$ then
$$
1 = \int g(-x) \sum_{\lambda\in\Lambda} f(x-\lambda) ~dx = 
 \sum_{\lambda\in\Lambda} \int g(-x) f(x-\lambda) ~dx,
$$
which, after the change of variable $y = -x+\lambda$, gives
$$
1 = \int f(-y) \sum_{\lambda\in\Lambda} g(y-\lambda) ~dy.
$$
This in turn implies,
since $\sum_{\lambda\in\Lambda}g(y-\lambda) \le 1$,
that $\sum_\lambda g(y-\lambda) = 1$ for
a.e.\ $y \in -\supp f$.

To complete the proof of the theorem, notice that if $f+\Lambda$ is a tiling
of $\RR^d$ and $a \in \RR^d$ is arbitrary then
both $f(x-a) + \Lambda$ and $g(x-a) + \Lambda$ are packings and
$f+\Lambda$ tiles $-\supp g(x-a) = -\supp g - a$.
We conclude that $g(x-a) + \Lambda$ tiles $-\supp f$, or
$g + \Lambda$ tiles $-\supp f - a$. Since $a\in\RR^d$ is arbitrary
we conclude that $g + \Lambda$ tiles $\RR^d$.
\Qed

\noindent
{\bf Example:}
Use Lemma \ref{lm:tiling-from-packing}
to prove that there is no measurable nonnegative function
$f$ that tiles with $\Lambda=\ZZ^d\setminus\Set{0}$ (or even $\ZZ^d$
minus a set of lower density $0$, such as a line).
Try to prove this otherwise.

\subsubsection{Failure of the lemma for non-translational tiling}
\label{sec:counterexample}

Suppose we study tiling where all rigid motions of the tile,
and not just translations, are allowed.
The analogue of the tiling set then is a set $\Lambda$ of rigid motions.
For $x\in\RR^d$ and $\lambda$ a rigid motion we denote by $\lambda(x)$ the action
of $\lambda$ on $x$.
The following theorem shows that our Lemma \ref{lm:tiling-from-packing}
is very particular to translations.
\begin{figure}[h] \centering
\input{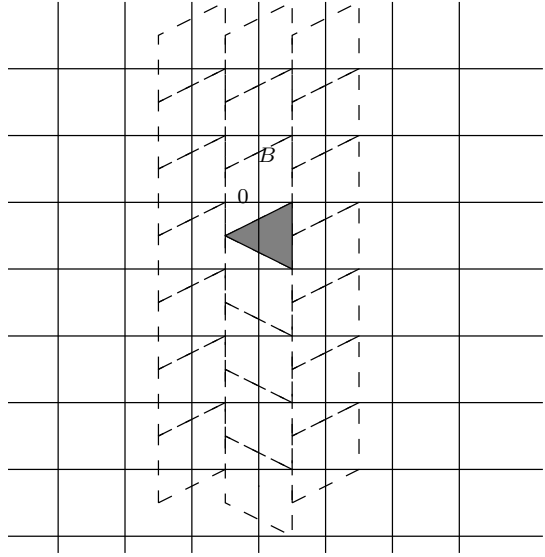}
\caption{Packing of set $B$, the parallelogram above
the shaded triangle, with motions $\Lambda$. The shaded triangle is not covered.}
\ \\ \hrule
\label{fig:counterexample}
\end{figure}

\begin{theorem}\label{counterexample}
There are two polygons $A$ and $B$ in $\RR^2$ of the same area
and a set of rigid motions $\Lambda$ such that both
collections
$\Set{\lambda(A):\ \lambda\in\Lambda}$ and
$\Set{\lambda(B):\ \lambda\in\Lambda}$ are packing
but only one of them is a tiling.
\end{theorem}

\noindent{\bf Proof.}
Take $A = (-1/2, 1/2)^2$
and $B$ to be the parallelogram with vertices
$(-1/2, -1/2)$, $(1/2, 0)$, $(1/2, 1)$ and $(-1/2, 1/2)$.
Take the set of rigid motions to be
the set of translations by $\ZZ^2$ modified as follows:
instead of translating by the elements $(0,k)$, $k<0$, we
first reflect the domain with respect to the $x$-axis and then
translate it by $(0,k)$. For the elements $(m,n)$ of $\ZZ^2$
where either $m\neq 0$ or $n\ge 0$ we just translate.

Since the reflection has no effect on $A$ the collection
$\Set{\lambda(A):\ \lambda\in\Lambda}$ clearly constitutes a tiling.
On the other hand the collection $\Set{\lambda(B):\ \lambda\in\Lambda}$
can be seen in Figure \ref{fig:counterexample} and is clearly not a tiling,
although it is a packing.
\Qed
\subsubsection{Deducing tiling from the condition on supports}

Assume that we have
\beql{the-usual-condition}
\supp \ft{\delta_\Lambda} \subseteq \Set{\ft f = 0} \cup \Set{0}
\eeq
for some non-zero $f\ge 0$ in $L^1$ and that $\Lambda$ is of bounded density.
Since $\ft f(0) = \int f > 0$ it follows that in some neighborhood $N$ of $0$
we have $(\supp\ft{\delta_\Lambda}) \cap N = \Set{0}$.
Hence the set
\beql{O-set}
O = \left( \supp \ft{\delta_\Lambda} \setminus \Set{0} \right)^c
\eeq
is open and
$$
\Set{\ft f \neq 0} \subseteq O.
$$

We shall need the following result.
\begin{theorem}\label{th:supp-implies-tiling}
Suppose that $0\le f \in L^1(\RR^d)$, $\int f =1$,
$\Lambda$ (of uniformly bounded density) is of density $1$,
and that \eqref{the-usual-condition} holds.
Suppose also that for the open set $O$ of \eqref{O-set} and for each
$\epsilon>0$ there exists $f_\epsilon \ge 0$ in $L^1(\RR^d)$ such that
$\ft{f_\epsilon}$ is in $C^\infty$, $\supp \ft{f_\epsilon} \subseteq O$ and
$$
\Lone{f - f_\epsilon} \le \epsilon.
$$
Then $f + \Lambda$ is a tiling.
\end{theorem}

\noindent
{\bf Proof.}
Suppose that $f_\epsilon$ is as in the Theorem.
First we show that $(\int f_\epsilon)^{-1}f_\epsilon + \Lambda$
is a tiling.
That is, we show that the convolution $f_\epsilon * \delta_\Lambda$
is a constant.
Let $\phi$ be $C_c^\infty$ function.
Then
$$
(f_\epsilon * \delta_\Lambda)(\phi) = 
\ft{f_\epsilon}\ft{\delta_\Lambda}(\ft\phi) =
\ft{\delta_\Lambda}(\ft\phi \ft{f_\epsilon}).
$$
But the function $\ft\psi = \ft\phi \ft{f_\epsilon}$ is a $C_c^\infty$
function whose support intersects $\supp\ft{\delta_\Lambda}$ only at $0$.
And, it is not hard to show, because $\Lambda$ has density $1$,
that $\ft{\delta_\Lambda}$ is equal to $\delta_0$ in a
neighborhood of $0$ (see \cite{kolountzakis:polygons}).
Hence 
$$
(f_\epsilon * \delta_\Lambda)(\phi) =
 \left(\ft\phi \ft{f_\epsilon}\right) (0) =
 \int \phi \int f_\epsilon,
$$
and, since this is true for an arbitrary $C_c^\infty$ function $\phi$,
we conclude that $f_\epsilon * \delta_\Lambda = \int f_\epsilon$,
as we had to show.

For any set $\Lambda$ of uniformly bounded density we have ($B$ is any ball
in $\RR^d$ and $g\in L^1(\RR^d)$)
$$
\int_B\Abs{\sum_{\lambda\in\Lambda} g(x-\lambda)} ~dx \le C_{B,\Lambda}
	\int_{\RR^d}\Abs{g},
$$
(See \cite{kolountzakis-lagarias} for a proof of this in dimension $1$, which holds
for any dimension.)
Applying this for $g=f-f_\epsilon$ we obtain that
$$
\sum_{\lambda\in\Lambda} f_\epsilon(x-\lambda) \to
 \sum_{\lambda\in\Lambda} f(x-\lambda), \ \ \ \mbox{in $L^1(B)$}.
$$
Since $B$ is arbitrary this
implies that $\sum_{\lambda\in\Lambda} f(x-\lambda) = 1$, a.e.\ in $\RR^d$.
\Qed

We write $\widetilde f(x) = \overline{f(-x)}$.

Let $\Omega \subset \RR^d$ be a bounded open set of measure $1$,
$\chi_\Omega$ its indicator function and
$f$ be such that $\ft{f} = \chi_\Omega * \widetilde{\chi_\Omega}$.
Then $\widetilde{f} = \Abs{\ft{\chi_\Omega}}^2 \ge 0$,
$\int f = 1$ by Parseval's theorem.
Clearly we have $\Set{\ft{f} \neq 0} = \Omega - \Omega$.

Write
$$
\Omega_\epsilon = \Set{x\in\Omega:\ {\rm dist\,}(x, \partial\Omega) >
\epsilon},
$$
and define $f_\epsilon$ by
$$
\ft{f_\epsilon} = \psi_\epsilon * \chi_{\Omega_\epsilon} *
 (\psi_\epsilon * \chi_{\Omega_\epsilon})^{\widetilde{\ }}
$$
(or $\widetilde{f_\epsilon} = \Abs{\ft{\psi_\epsilon}}^2
 \Abs{\ft{\chi_{\Omega_\epsilon}}}^2$),
where $\psi_\epsilon$ is a smooth, positive-definite
approximate identity supported in $B_{\epsilon/2}(0)$.

One can easily prove the following proposition.

\noindent
{\em If $g_n \to g$ in $L^2$ then $\Abs{g_n}^2 \to \Abs{g}^2$ in $L^1$.}

(For the proof just notice the identity
$$
\Abs{g}^2 - \Abs{g_n}^2 = \Abs{g - g_n}^2 +
  2\cdot{\rm Re} \left(\overline{g_n}(g-g_n)\right),
$$
integrate and use the triangle and Cauchy-Schwartz inequalities.)

Since $\psi_\epsilon*\chi_{\Omega_\epsilon} \to \chi_\Omega$ in $L^2$
(dominated convergence) we have (Parseval)
that $\ft{\psi_\epsilon}\ft{\chi_{\Omega_\epsilon}} \to \ft{\chi_\Omega}$
in $L^2$ and, using the proposition above, that
$\Abs{\ft{\psi_\epsilon}}^2 \Abs{\ft{\chi_{\Omega_\epsilon}}}^2
  \to \Abs{\ft{\chi_\Omega}}^2$ in $L^1$,
which means that $f_\epsilon \to f$ in $L^1$.

We also have that
$$
\supp{\ft{f_\epsilon}} \subseteq
 \overline{\Omega_{\epsilon/2}} - \overline{\Omega_{\epsilon/2}}
\subseteq \Omega - \Omega = \Set{\ft f \neq 0}.
$$

The assumptions of Theorem \ref{th:supp-implies-tiling} are therefore
satisfied.
Combining Theorems \ref{th:tiling-implies-supp} and
\ref{th:supp-implies-tiling} with the above observations we
obtain the following characterization of tiling
by the function $\Abs{\ft{\chi_\Omega}}^2$.
The special form of this function allows us to drop
any conditions, that are otherwise needed, regarding the order (how many derivatives
it involves) of the
tempered distribution $\ft{\delta_\Lambda}$.
\begin{theorem}\label{th:supp-equiv-tiling}
Let $\Omega$ be a bounded open set, $\Lambda$ a discrete set in $\RR^d$,
and $\delta_\Lambda = \sum_{\lambda\in\Lambda} \delta_\lambda$.
Then $\Abs{\ft{\chi_\Omega}}^2 + \Lambda$ is a tiling if and only if
$\Lambda$ has uniformly bounded density and
$$
(\Omega - \Omega) \cap \supp\ft{\delta_\Lambda} = \Set{0}.
$$
\end{theorem}

\noindent
{\bf Proof of Theorem \ref{th:spectra-of-the-cube}.}
By a simple calculation we get
\begin{eqnarray*}
{\cal Z}(\ft{\chi_Q}) &=& \Set{\xi\in\RR^d:\ \mbox{some $\xi_j$
  is a non-zero integer}}\\
 &\subseteq& (2Q)^c.
\end{eqnarray*}

Suppose first that $Q+\Lambda=\RR^d$.
From Theorem \ref{th:tiling-implies-supp} it follows that
\begin{eqnarray*}
\supp\ft{\delta_\Lambda} &\subseteq& \Set{0}\cup {\cal Z}(\ft{\chi_Q})\\
 &\subset& \Set{0}\cup (Q-Q)^c
\end{eqnarray*}
and from Theorem \ref{th:supp-equiv-tiling} we deduce that $\Lambda$
is a spectrum of $Q$.

Conversely assume that $\Lambda$ is a spectrum of $Q$, so that
$\Abs{\ft{\chi_Q}}^2+\Lambda=\RR^d$.
It follows that $(Q-Q)\cap(\Lambda-\Lambda) = \Set{0}$ as we
have $\Abs{\ft{\chi_Q}}^2(0)=1$ and
$\Abs{\ft{\chi_Q}}^2>0$ on $Q-Q$.
But this means that we have a packing $Q+\Lambda \le \RR^d$.
However, $\Lambda$ is a tiling set, because it is a spectrum, and
there is another object that tiles with $\Lambda$, namely
$\Abs{\ft{\chi_Q}}^2$, and this object has the same integral
as $\chi_Q$ (that is, $1$).
It follows from Lemma \ref{lm:tiling-from-packing} that
$Q+\Lambda = \RR^d$ is also a tiling, as we had to prove.
\Qed

\subsection{A proof that the disk is not spectral, which just makes it}
\label{sec:bessel}

Here we present a proof of why the disk $D=\Set{\Abs{x}<{1\over\sqrt\pi}}$ in
the plane is not a spectral domain.
The radius is taken equal to $1/\sqrt\pi$ to
make the disk have area $1$, as we usually do in this survey.

The proof is simple but relies on two not-so-easy facts.
\begin{enumerate}
\item\label{en1}
The first is the upper bound ${\pi\over\sqrt{12}}$, due to Thue,
on the density of any
packing of the plane with copies of the same disk
(see, for example, \cite[Ch.\ 3]{pach-agarwal}).
\item\label{en2}
The second is that the first zero of the Fourier Transform of the indicator
function of $D$ is at distance approximately 1.08098 from the origin. This
may either be looked up in tables of the Bessel function $J_1$ (which, up to
scaling, is the Fourier Transform of the indicator function of $D$ restricted on a line) or
may be computed in a straightforward way using a computer.
(The Fourier Transform of the unit-area disk,
defined by $\ft{\chi_D}(\xi) = \int_D \exp(-2\pi i \xi x)~dx$,
is equal to a constant times
$J_1(2\sqrt\pi \Abs{\xi})$ and the first zero of $J_1$ is at $3.832\cdots$.)
\end{enumerate}
Fuglede \cite{fuglede} was the first to suggest that the disk is not spectral,
but the argument was unclear.
The situation has since been
clarified in the papers of Iosevich, Katz and Pedersen \cite{iosevich-katz-pedersen},
who proved that the ball in any dimension is not spectral, and of
Iosevich, Katz and Tao \cite{iosevich-katz-tao}, in which a much more general result is
proved: every smooth convex hypersurface cannot have an interior which is a
spectral domain.
It was also shown by Fuglede \cite{fuglede:ball} (for the Euclidean ball in $\RR^d$)
and by Iosevich and Rudnev \cite{iosevich-rudnev} (for any smooth convex body in $\RR^d$,
for $d \neq 1 \bmod 4$)
that there can only be a finite number of orthogonal exponentials in the corresponding $L^2$ spaces.

The method shown in this section is still
interesting because of its simplicity and, perhaps, entertaining as the fact that it works appears to be an accident.

\begin{figure}[h] \centering
\input{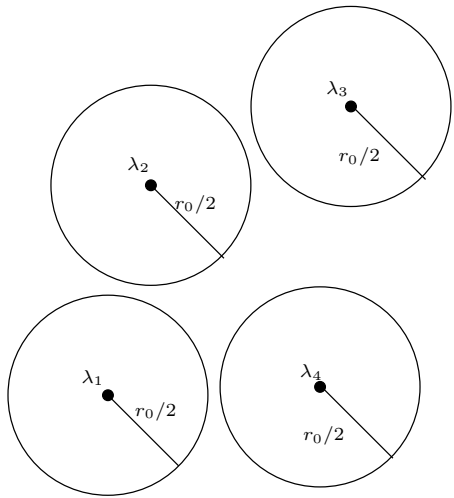}
\caption{The packing by disks of radius $r_0/2$ centered at the points of the spectrum.
Thue's result means that the area outside the disks has density $\ge 1-\pi/\sqrt{12}$}
\ \\ \hrule
\label{fig:packing}
\end{figure}

The Fourier Transform of $D$ is radial, as is
the function itself, hence the set of zeros of the Fourier Transform is a set of circles
centered at the origin. Let $r_0$ be the radius of the smallest such circle. By a simple
numerical calculation we locate $r_0 = 1.08098\cdots$. Suppose now that the disk is
spectral with spectrum $\Lambda$. 
Since $\Lambda-\Lambda \subseteq \Set{\ft{\chi_D} = 0} \cup \Set{0}$
it follows that $\Abs{\lambda - \mu} \ge r_0$ for any $\lambda, \mu \in \Lambda$,
$\lambda \neq \mu$, and hence, if we
center a copy of a disk of radius $r_0/2$, call it $D_1$,  at each point of
$\Lambda$, we have a packing of the plane with congruent disks (see Figure \ref{fig:packing}).
The density of
such a packing is at most $\pi/\sqrt{12}$, by Fact 1 above.

Since the integral of the {\em power spectrum} $\Abs{\ft{\chi_D}}^2$ of $\chi_D$ is $1$
(Parseval), and the power spectrum tiles with $\Lambda$
it follows that the density of $\Lambda$ is
equal to $1$ as well, hence the density of the packing $D_1 + \Lambda$
is equal to the area of $D_1$, which is $\pi r_0^2/4$. So we have the
inequality
$$
\pi {r_0^2 \over 4} \le {\pi \over \sqrt{12}},
$$
which implies
$$
r_0 \le {2 \over (12)^{1/4}} = 1.0745699\cdots,
$$
which is in contradiction with Fact \ref{en2} above which states that $r_0$ is
approximately 1.08098.

\subsection{More results on the Fuglede Conjecture}

\subsubsection{Convex domains}

The convex bodies which tile space have long been known \cite{venkov,mcmullen}
to be precisely the polytopes which are symmetric, have symmetric co-dimension one
facets and their co-dimension two facets each have a ``belt'' which consists
of four or six facets (the belt of a facet is the collection of all facets
of the polytope which are translates of the given facet).
It is also known \cite{mcmullen} that whenever a convex body $\Omega$ tiles space by translation
it can also tile by lattice translation.
It follows from Theorem \ref{th:fuglede-lattice} that convex bodies which tile
are also spectral, and possess a lattice spectrum (the dual lattice of
their translation lattice).

Our knowledge is much less complete for convex bodies which are spectral.
In particular we do not know yet that spectral convex bodies are also tiles,
but we are getting there.
Most of the results described in this section are in the general direction
of showing that well known facts which hold for convex tiles are also
true of convex spectral bodies.

In \cite{iosevich-katz-tao} it was proved that smooth convex bodies
cannot be spectral, a fact which is clearly true of convex bodies which tile,
even if one has not heard of the Venkov-McMullen theorem.
\begin{theorem}\label{th:smooth-not-spectral}
{\bf(Iosevich, Katz and Tao, 1999)}\\
Suppose that  $\Omega$ is a symmetric convex body in $\RR^d$, $d\ge 2$.
If the boundary of $\Omega$  is smooth, then it  does not admit a spectrum.
The same conclusion holds in $\RR^2$ if the boundary is piecewise smooth,
and has at least one point of non-vanishing Gaussian curvature.
\end{theorem}
The starting point of the proof is the fact that the zero set
$$
{\cal Z} = \Set{\ft{\chi_\Omega} = 0}
$$
is known, asymptotically, to an ever-higher degree of accuracy.
For example, it is a well known fact (see e.g.\ \cite{iosevich-katz-tao})
that if $\xi$ is a zero of $\ft{\chi_\Omega}$ and $\xi\to\infty$
such that $\xi$ remains inside a cone
$$
{\cal C} = \Set{\xi:\ {\inner{\xi}{u} \over \Abs{\xi}} > 1-\epsilon},
$$
where $u \in S^{d-1}$ is the unit outward normal vector at some point
$x \in \partial\Omega$ of positive curvature and $\epsilon>0$ is sufficiently
small, then
$$
\Norm{\xi}_{\Omega^o} = \left({\pi\over 2} + {d\pi \over 4}\right) + k\pi + o(1),
        \ \ \ (\xi\to\infty),
$$
where $\Omega^o$ is the dual body (which is also smooth), $d$ is the dimension
and $k$ is an integer.
One then uses the fact that if $\Lambda$ is a spectrum then
$\Lambda-\Lambda \subseteq {\cal Z}$ in order to reach a contradiction.

It turns out \cite{iosevich-rudnev} that for smooth
convex bodies with nowhere vanishing
Gaussian curvature (such as the Euclidean ball) much more is
true than the fact that there is no complete orthogonal set
of exponentials for their $L^2$ space.
\begin{theorem}\label{th:finitely-many-exponentials}
{\bf(Iosevich and Rudnev, 2002)}\\
Suppose that  $\Omega$ is a smooth symmetric convex body in $\RR^d$, $d\ge 2$,
with nowhere vanishing Gaussian curvature. If $d \neq 1 \bmod 4$ then
any set of orthogonal exponentials in $L^2(\Omega)$ is finite.
If $d=1 \bmod 4$ such a set may be infinite only if it is a subset
of a one-dimensional lattice.
\end{theorem}
This has also been proved for the ball in any dimension by Fuglede \cite{fuglede:ball}.

Finally, in dimension $d=2$ the Fuglede Conjecture may be considered
settled for convex bodies \cite{iosevich-katz-tao:plane}.
\begin{theorem}\label{th:fuglede-plane}
{\bf(Iosevich, Katz and Tao, 2002)}\\
The only convex domains in $\RR^2$ which are spectral
are the parallelograms and the symmetric hexagons
(these are the only convex tiles as well).
\end{theorem}

\subsubsection{Polytopes with unbalanced facets}

Suppose that $\Omega$ is a polytope, not necessarily convex, that tiles
space by translation. Suppose also that $u$ is one of its face normals
and let $F_1^+, \ldots, F_k^+$ be all its facets with outward normal
in the direction of $u$ and let $F_1^-,\ldots,F_l^-$ be the facets
with outward normal in the direction of $-u$.
One can easily see that we must have
$$
\Abs{F_1^+}+\cdots+\Abs{F_k^+} = \Abs{F_1^-}+\cdots+\Abs{F_l^-}.
$$
The reason is that in any tiling by translates of $\Omega$
the facets $F_j^+$ can only be ``countered'' by translates of the facets
$F_j^-$.
Applying this for a large region in space one deduces that the total area
of the plus-facets must equal that of the minus-facets. 

The following result \cite{kolountzakis-papadimitrakis:polytopes}
claims that spectral polytopes have the same property.
\begin{theorem}\label{th:unbalanced-polytopes}
{\bf(Kolountzakis and Papadimitrakis, 2000)}\\
If $\Omega$ is a polytope in $\RR^d$ which, for some direction $u$ normal
to a facet, has more area with outward normal $u$ than it has
with outward normal $-u$, then $\Omega$ is not spectral.
Clearly it can also not be a tile.
\end{theorem}
\begin{figure}[h] \centering
\input{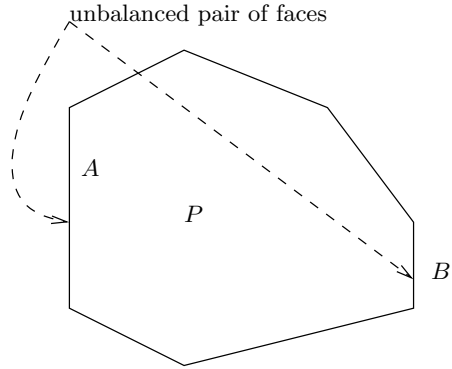}
\caption{A polytope $P$ with many directions of unbalanced faces. The two
facets shown are the only ones perpendicular to their normal, yet
there is more face measure looking left than looking right.
Such a polytope can neither tile by translation nor be spectral.}
\ \\ \hrule
\label{fig:unbalanced}
\end{figure}
We do not present the proof of this result here. However, the following
toy-case is rather instructive.
Suppose that we have a polytope $\Omega$ which has precisely two facets
$A$ and $B$ (see the example in Figure \ref{fig:unbalanced})
with normals parallel to a certain $u \in S^{d-1}$.
Assume that facet $A$ has outward normal $u$ and facet $B$ has $-u$, and
that the area of $A$ is not equal to that of $B$.

We claim that in any semi-infinite tube whose axis is the line
$\RR u$ and {\em any} bounded domain as base there are only finitely
many points of any spectrum. This is impossible as for any spectrum there
is a number $R$ such that in any ball of radius $R$ we can find some
point of the spectrum.
To show the above claim it is enough to show that any such tube is eventually
(that is, near infinity) free from zeros of $\ft{\chi_\Omega}$, or, what
amounts to the same thing, free from zeros of
$$
\ft{\nabla_u\chi_\Omega}(\xi) = 2\pi i \inner{\xi}{u} \ft{\chi_\Omega}(\xi).
$$
Observe now that $\nabla_u\chi_\Omega$ is a measure supported
on the facets of the polytope, which is a constant function on
every facet, a constant which depends on the angle the facet
is forming with $u$.

Look then at what happens to the Fourier Transform $\ft{\nabla_u\chi_\Omega}$
along the line $\RR u$. Along that line the values of the Fourier Transform
that we are reading are just the values of the {\em one-dimensional}
Fourier Transform of the {\em projection} of the measure
$\nabla_u\chi_\Omega$ on the line $\RR u$.
This is the measure $\mu$ defined by
$$
\mu(E) = \nabla_u\chi_\Omega(E + u^\perp),\ \ \ (E \subseteq \RR),
$$
and it is clear that $\mu$ has a continuous part coming from all the
facets which are non-orthogonal to $u$ and also contains the two point
masses $\Abs{A}\delta_a$ and $-\Abs{B}\delta_b$, where $a, b \in \RR$
are the points on $\RR u$ where the facets $A$ and $B$ project.
By the Riemann-Lebesgue lemma the contribution to $\ft\mu$ of the
continuous part of $\mu$ fades to $0$ as we tend to $\infty$ and
it is the Fourier Transform of the atomic part that dominates $\ft\mu$, namely
(as $t\to\infty$)
$$
\ft\mu(t) \sim \Abs{A}e^{2\pi i \inner{a}{t}}-\Abs{B}e^{2\pi i \inner{b}{t}}
$$
whose absolute value is $\ge\Abs{\Abs{A}-\Abs{B}}$.
So, for large $t$, there are no zeros on the line, and with a little more care,
we can show that the same (albeit farther away) is true in any tube around this
line.

\subsubsection{Dimension $1$}

Even in dimension $1$ the Fuglede Conjecture appears to be rather hard.
The number-theoretic aspect of the problem is seen more clearly
here, especially if one looks just at sets of the type
$$
\Omega = A + (0,1),\ \ \ \mbox{($A$ a finite subset of $\ZZ$)}.
$$
The conjecture is still open for this class of sets.

The following are interesting partial results.
\begin{enumerate}
\item
\L aba \cite{laba:two} showed that whenever $\Abs{A}=2$ the conjecture
is true.
\item
This is improved to $\Abs{A}=3$ by \L aba in \cite{laba:multi}.
In the same paper it is also shown that if $\Abs{A}$ has at most
two prime factors then if $\Omega$ is a tile it is also spectral.
\item
\L aba also shows in \cite{laba:multi} that
if $\Abs{A} > 3/2(\max A - \min A)$ then the $\Omega$ is a tile
if and only if it is spectral.
This is generalized by Kolountzakis and \L aba \cite{kolountzakis-laba}
to any set $\Omega$ of measure $1$ which is a subset of $(0,3/2-\epsilon)$,
for some $\epsilon>0$.
In fact what is really shown in \cite{kolountzakis-laba} is that such ``tight'' domains
can only be spectral or tiles if they tile by the lattice $\ZZ$.
\end{enumerate}

\end{document}